\documentclass[11pt]{article}
\usepackage{amsmath,amssymb,amsthm}
\usepackage{authblk}
\usepackage[utf8]{inputenc}
\usepackage{fontenc}
\usepackage{mathtools, nccmath, etoolbox}
\usepackage[a4paper]{geometry}
\usepackage[hidelinks]{hyperref} 
\usepackage{graphicx}
\usepackage[numbers]{natbib} 
\usepackage{mathdots}
\usepackage[monochrome]{xcolor} 
\usepackage{comment}
\usepackage{tikz}	
\usepackage{pgf}
\usepackage{pgfplots}
\pgfplotsset{compat=1.18}
\usepackage{caption}
\usepackage[font=footnotesize,labelformat=simple]{subcaption}
\usepackage{hyperref}

\theoremstyle{plain}

\theoremstyle{definition}

\theoremstyle{remark}
\newtheorem*{remark*}{Remark}

\begin{document}

\title{Discrete Adjoint Method for Variational Integration of Constrained ODEs and its application to Optimal Control of Geometrically Exact Beam Dynamics}

\author[1]{Matthias Schubert \thanks{matthias.schubert@fau.de}}
\author[1]{Rodrigo T. Sato Mart{\'\i}n de Almagro \thanks{rodrigo.t.sato@fau.de}}
\author[2,3]{Karin Nachbagauer \thanks{karin.nachbagauer@fh-wels.at}}
\author[4]{Sina Ober-Blöbaum \thanks{sinaober@math.uni-paderborn.de}}
\author[1]{Sigrid Leyendecker \thanks{sigrid.leyendecker@fau.de}}

\affil[1]{Friedrich-Alexander-Universität Erlangen-Nürnberg, Institute of Applied Dynamics, Immerwahrstrasse 1, 91058 Erlangen \vspace*{0.2cm}}
\affil[2]{Faculty of Engineering and Environmental Sciences, University of Applied Sciences Upper Austria, Stelzhamerstraße 23, 4600 Wels, Austria \vspace*{0.2cm}}
\affil[3]{Institute for Advanced Study, Technical University of Munich,  Lichtenbergstraße 2a, 85748 Garching, Germany \vspace*{0.2cm}}
\affil[4]{Universität Paderborn, Warburger Straße 100, 33098 Paderborn\vspace*{0.2cm}}
\date{}

\maketitle

\begin{abstract}
Direct methods for the simulation of optimal control problems apply a specific discretization to the dynamics of the problem, and the discrete adjoint method is suitable to calculate \textcolor{blue}{corresponding conditions to approximate an optimal solution. While the benefits of structure preserving or geometric methods have been known for decades, their exploration in the context of optimal control problems is a relatively recent field of research.} In this work, the discrete adjoint method is derived for variational integrators \textcolor{blue}{yielding structure preserving approximations of the dynamics firstly} in the ODE case and \textcolor{blue}{secondly} for the case in which the dynamics is subject to holonomic constraints. The convergence rates are illustrated by numerical examples. \textcolor{blue}{Thirdly,} the discrete adjoint method is applied to geometrically exact beam dynamics, represented by a holonomically constrained PDE.
\end{abstract}

\bigskip
\noindent This work was published in Multibody System Dynamics on the 5th of September 2023. \href{https://doi.org/10.1007/s11044-023-09934-4}{https://doi.org/10.1007/s11044-023-09934-4}




\section{Introduction}
\label{sect:Intro}
\textcolor{blue}{There are two alternative ways to handle an optimal control problem numerically. The so-called \emph{indirect methods} first derive the necessary conditions for optimality in the continuous-time setting by applying \textrm{\textsc{Pontryagin}}'s maximum principle and then discretizing the resulting equations. In contrast, \emph{direct methods} first discretize the continuous problem, turning it into a finite dimensional one, and then apply a discrete version of \textrm{\textsc{Pontryagin}}'s maximum principle. In both cases, one is led to the augmentation of the original objective with the different constraints enforced by \textrm{\textsc{Lagrange}} multipliers. The \textrm{\textsc{Lagrange}} multipliers enforcing the plant (the dynamic equations) of the problem are commonly called adjoint or co-state variables. In the multibody systems literature it is common to refer to this as the adjoint method, and in particular, the discrete adjoint method when considered as a direct method.} In this contribution, we \textcolor{blue}{apply} the discrete adjoint method to optimal control problems with variational integrators approximating the dynamics.\\
In general for direct approaches, the discretization of the ODE governing the dynamics results in a specific discretization of the adjoint variables especially for symplectic methods as e.g. variational integrators \cite{Campos2015,Marsden2011,Bonnans2006}. Variational and thus symplectic numerical methods are worthy of consideration as they can benefit the solution of boundary value problems \cite{Offen2018}. For the optimal control of constrained ODEs, discretizations with conservation properties are of interest as well \cite{Betsch2021,Leyendecker2011}.\\
The optimal control of mechanical PDEs, such as string and beam dynamics is an active field of research \cite{Betsch2022,Lismonde2019,Seifried2014}.
The discrete adjoint method has been used for the optimization of flexible multibody systems \cite{Callejo2019} as well as for parameter identification in rigid body dynamics \cite{Lauss2018,Lauss2017}.
The discrete adjoint method for variational integrators with holonomic constraints is discussed in \cite{Ebrahimi2019}.
The discrete adjoint method is derived for a specific discretization of dynamics and this matches the chosen integrator. Therefore, it suggests itself to be applied to integrators that are structure preserving \cite{Sanz2016}.\\
In this work we briefly summarize variational integrators and then show how to derive the discrete adjoint equations for this class of integrators. The basic principles, the derivation of boundary conditions and the discretization of forces are explained. The discrete adjoint method is then extended to variational integrators for holonomically constrained ODEs. The convergence behavior of both methods is investigated with the example of a mathematical pendulum. Finally, the method is applied to the constrained PDE case of geometrically exact beam dynamics.

\section{Discrete Adjoint Method for Variational Integrators}\label{sect:DAMVIs}
\subsection{Variational Integrators}\label{sect:VI}%
This section illustrates the derivation of the equations of motion for forced systems via variational principles in the continuous and discrete setting \cite{Marsden2001,Marsden2011}. These equations have to be fulfilled as constraints for the optimal control problem.\\
\textcolor{blue}{Consider a Lagrangian mechanical system whose configuration space is the $n$-dimensional smooth manifold $Q$. The motion of our system is represented by a curve $q: [0, T] \to \mathcal{Q}$, $t \mapsto q(t)$. We denote the velocity of the configuration at time $t$ by $\dot{q}(t) \in \mathcal{T}_{q(t)}\mathcal{Q}$. The Lagrangian is a function defined on the tangent bundle of $\mathcal{Q}$, $\mathcal{TQ}$, $\mathcal{L}: \mathcal{TQ} \to \mathbb{R}$. It} usually represents the difference of kinetic and potential energy. \textcolor{blue}{An external Lagrangian control force is a map $f_\mathcal{L}: \mathcal{TQ} \times \mathcal{U} \to \mathcal{T}^*\mathcal{Q}$ where $\mathcal{U} \subseteq \mathbb{R}^l$, $l \leq n$, is the space of admissible controls. A control is thus a curve $u: [0,T] \to \mathcal{U}$. The total virtual work of such a system vanishes}
\begin{equation}
  \delta \int_0^T \mathcal{L}\big(q(t), \dot{q}(t) \big)~dt + \int_0^T f_\mathcal{L} \big(q(t), \dot{q}(t), u(t)\big) \, \delta q\textcolor{blue}{(t)} ~dt =0,\quad \forall \textcolor{blue}{\delta q(t)}
\end{equation}
\textcolor{blue}{This is the \textrm{\textsc{Lagrange}}-\textrm{\textsc{d'Alembert}} principle (with controls), which states that the total virtual work evaluated over a physical trajectory of the system $q$ (and a control $u$) vanishes for all variations $\delta q(t)$ with fixed end-points $\delta q(0) = \delta q(T) = 0$}. This leads to the equations of motion, the forced \textrm{\textsc{Euler}}-\textrm{\textsc{Lagrange}} equations:
\begin{equation}
    -\frac{d}{dt} \frac{\partial \mathcal{L}(q, \dot{q})}{\partial \dot{q}} + \frac{\partial \mathcal{L}(q, \dot{q})}{\partial q} + f_\mathcal{L}(q, \dot{q}, u) = 0.
\end{equation}
This principle is an extension of \textrm{\textsc{Hamilton's}} principle to include non-\-con\-ser\-va\-ti\-ve forces such as control or dissipative forces. 
A forced variational integrator is derived via the approximation of the action and the virtual work of non-conservative forces and subsequent variation in the discrete setting \cite{Marsden2001,Hairer2010,Jordan1964}. The time interval $[0,T]$ is discretized by $N$ time nodes, we consider a discrete configuration path $\{ q_n \}_{n=0}^N$ with $q_n \approx q(t_n)$ with linear approximation of $q(t)$ in $[t_n, t_{n+1}]$. The approximation of the action integral via the discrete Lagrangian $L_d$ and the approximation of the virtual work of non-conservative forces via the left and right side discrete forces $f_d^-$ and $f_d^+$ is considered. The input variable is approximated as $u_n \approx u(t_n)$. In each time interval $[t_n, t_{n+1}]$, the control path $u_d = \{ u_n\}_{n=0}^{N-1}$ is approximated constant. 
\begin{align}
    \int_{t_n}^{t_{n+1}} \mathcal{L}\textcolor{blue}{(q(t),\dot{q}(t))} ~dt &\approx L_d(q_n, q_{n+1})\\
    \int_{t_n}^{t_{n+1}} f_\mathcal{L}\textcolor{blue}{(q(t),\dot{q}(t),u(t)) ~\delta q(t)} ~dt &\approx f_d^-(q_n, q_{n+1}, u_n) ~\delta q_n\nonumber\\
    &+ f_d^+(q_n, q_{n+1}, u_n) ~\delta q_{n+1} 
\end{align}
The \textcolor{blue}{discrete total} virtual work \textcolor{blue}{vanishes}:
\begin{equation}\label{eqn:LAprinc}
    \sum_{n=0}^{N-1} \textcolor{blue}{\left[ \delta  L_d(q_n, q_{n+1}) + f_d^-(q_n, q_{n+1}, u_n) ~\delta q_n + f_d^+(q_n, q_{n+1}, u_n) ~\delta q_{n+1} \right]}=0, \; \forall \delta q_n
\end{equation}
with $\delta q_0 = \delta q_N = 0$. The discrete \textrm{\textsc{Lagrange}}-\textrm{\textsc{d'Alembert}} principle leads to the discrete, forced \textrm{\textsc{Euler}}-\textrm{\textsc{Lagrange}} equations, which are derived via discrete variation and subsequent rearrangement of terms for fixed boundary conditions. The slot derivatives $D_k$ denote derivatives with respect to the $k$-th argument.
\begin{equation}\label{eqn:DEL}
    D_1 L_d(q_n, q_{n+1}) + D_2 L_d(q_{n-1}, q_n) + f_d^-(q_n, q_{n+1}, u_n) + f_{d}^+(q_{n-1}, q_n, u_{n-1}) = 0,
\end{equation}
for $n=1,~...,~N-1$. This equation takes two positions at the current and the previous time node and defines the relation with the next one. Given $q_{n-1}, q_n, u_{n-1}$ and $u_n$, \textcolor{blue}{this equation determines a unique $q_{n+1}$ provided the discrete Lagrangian is regular, i.e. the matrix $D_1 D_2 L_d = D_2 D_1 L_d$ is regular.}

The initial conditions are usually defined on $\mathcal{TQ}$ as position and velocity \textcolor{blue}{or on $\mathcal{T}^* \mathcal{Q}$ as position and momentum, but} not on $\mathcal{Q} \times \mathcal{Q}$ as two positions at different points in time. To initialize this time stepping scheme, \textcolor{blue}{both a continuous and discrete version of the \textrm{\textsc{Legendre}} transformation are needed}.

\textcolor{blue}{The continuous \textrm{\textsc{Legendre}} transformation, $\mathbb{F}\mathcal{L}: \mathcal{TQ} \to \mathcal{T^*Q}$, $(q,\dot{q}) \mapsto ( q, p = D_2 L(q,\dot{q}) )$ connects the Lagrangian and the Hamiltonian formulations of dynamics. It allows us to compute an initial momentum $p^0$ from an initial configuration and velocity, $(q^0, \dot{q}^0)$. In the discrete setting, the (forced) discrete \textrm{\textsc{Legendre}} transformation defines \textcolor{blue}{two distinct maps} from the discrete state space to the cotangent bundle, $\textcolor{blue}{\mathbb{F}^{\pm}}L_d: \mathcal{Q} \times \mathcal{Q} \times \mathcal{U} \to \mathcal{T}^*\mathcal{Q}$, defined by}
\begin{subequations}\label{eqn:Leg}
    \begin{align}
    \mathbb{F}^- L_d: (q_n, q_{n+1}, u_n) \mapsto &(q_n, p_n^-) \\
    &= \Big(q_n, -D_1 L_d(q_n, q_{n+1}) - f_d^-(q_n, q_{n+1}, u_n) \Big) \nonumber\\
    \mathbb{F}^+ L_d: (q_n, q_{n+1}, u_n) \mapsto &(q_{n+1}, p_{n+1}^+)\\
    &= \Big(q_{n+1}, D_2 L_d(q_n, q_{n+1}) + f_d^+(q_n, q_{n+1}, u_n)\Big),\nonumber
    \end{align}
\end{subequations}
with the left and right side discrete momenta $p_n^-$ and $p_n^+$. \textcolor{blue}{These allow us to interpret the discrete \textrm{\textsc{Euler}}-\textrm{\textsc{Lagrange}} equations (\ref{eqn:DEL}) as a matching of momenta $p_n^- = p_n^+$ for $n=1,~...,~N-1$.}

\textcolor{blue}{In order to initialize the algorithm, given a configuration $q^0$, a velocity $\dot{q}^0$ and an initial control $u_0$, the relation
\begin{equation}\label{eqn:LegIni}
    D_2 L(q^0, \dot{q}^0) = p^0 = -D_1 L_d (q^0, q_1) -f_d^-(q^0, q_1, u_0)
\end{equation}
determines $q_1$}.

\subsection{Derivation of the Discrete Adjoint Method for Variational Integrators}\label{sec:DAMVI}%

Similar to the discrete variational principle in Section \ref{sect:VI}, now the discrete adjoint method for variational integrators in \eqref{eqn:DEL} is derived via a discrete variational principle and the structure and the resulting numerical method for the adjoint equations are illustrated.\\ 
Here, we concentrate on a discrete objective $J_d$ containing a \textcolor{blue}{quadratic} \textrm{\textsc{Mayer}} term,
\begin{equation}\label{eqn:Mayerterm}
\textcolor{blue}{J_M(q_N,p_N)} = \textcolor{blue}{\frac{1}{2} (q_N-q^N)^T S_q (q_N-q^N) + \frac{1}{2}  (p_N-p^N)^T S_p (p_N-p^N)}
\end{equation}
\textcolor{blue}{where $S_q$ and $S_p$ are positive semidefinite matrices. The \textrm{\textsc{Mayer}} term is used to relax the enforcement of the end state conditions, $(q^N, p^N)$, introducing weights for the reaching of the configuration and the momentum at the last time step $N$.}

\textcolor{blue}{The discrete adjoint method is derived by augmenting the objective with the variational integrator (\ref{eqn:DEL}) and (\ref{eqn:LegIni}) as constraints and by taking variations of the augmented objective \cite{Marsden2011}.}

The resulting nonlinear constrained optimization problem reads
\textcolor{blue}{
\begin{subequations}\label{eqn:OCP}
    \begin{equation}
        \underset{\textcolor{blue}{q_d,} u_d}{\min}\, J_d(q_d, u_d) = J_M\textcolor{blue}{(q_N,p_N)} + \sum_{n=0}^{N-1} \textcolor{blue}{\frac{1}{2}} u_n^T R u_n
    \end{equation}
subject to:
    \begin{align}
        \hphantom{XXXXX}
        q_0 &= q^0,\\
        p^0 &= - D_1 L_d(q^0, q_1) - f_{d}^{-}(q^0, q_{1}, u_0),\\
        0 &= D_1 L_d(q_n, q_{n+1}) + f_{d}^{-}(q_n, q_{n+1}, u_n)\\
          &+ D_2 L_d(q_{n-1}, q_{n}) + f_{d}^{+}(q_{n-1}, q_{n}, u_{n-1}),\quad\text{for}~n=1,~...,~N-1 \nonumber\\
        p_N &= D_2 L_d(q_{N-1}, q_N) + f_{d}^{+}(q_{N-1}, q_{N}, u_{N-1}),\label{eqn:LegTran}
    \end{align}
\end{subequations}}
The quantities $p^0$ and $q^0$ are prescribed initial conditions at the initial time node. \textcolor{red}{The objective also includes a {\textrm{\textsc{Lagrange}}} term which is quadratic in the control and $R$ is a positive-definite weight matrix.} \textcolor{blue}{Equation \eqref{eqn:LegTran} defining $p_N$ corresponds to the discrete \textrm{\textsc{Legendre}} transformation $\mathbb{F}^+L_d(q_{N-1},q_{N},u_{N-1})$.}\\
\\
REMARK~1: \textcolor{blue}{The dependence on $q_{N-1}$ and $q_N$ of the momentum term (\ref{eqn:LegTran}) of the \textrm{\textsc{Mayer}} term makes it more prone to produce larger contributions than the configuration term. This can make the optimization process unstable and possibly not convergent. In order to improve this, an iterative approach may be used where the end momentum of the $(i)$-th iteration, $p_N^{(i)}$, is used to inform the choice of a modified desired end momentum, $\tilde{p}^N$, such that
\begin{equation*}
 \Vert p_N^{(i)}-\tilde{p}^N(p_N^{(i)},p^N) \Vert \leq \Vert p_N^{(i)}- p^N \Vert
\end{equation*}
with $\tilde{p}^N(p^N,p^N) = p^N$. The procedure can be initialized by considering a first iteration with $S_p = 0$, and ended once $\Vert p_N^{(i)}- p^N \Vert$ is sufficiently small to allow us to substitute $\tilde{p}^N$ by $p^N$ in a final iteration.} \\
\\
The objective $J_d$ is augmented to $\tilde{J}_d$ by the initial conditions and the discrete \textrm{\textsc{Euler}}-\textrm{\textsc{Lagrange}} equations via adjoint variables $\lambda_n \approx \lambda(t_n)$ with the discrete adjoint path $\lambda_d = \{\lambda_n\}_{n=0}^{N-1}$. The indices are chosen such that $\lambda_n$ pairs with the corresponding momenta $p^{\pm}_n$.
\begin{align}
        \tilde{J}_d(q_d, u_d, \lambda_d) &= J_M\textcolor{blue}{(q_N,p_N(q_{N-1},q_N,u_{N-1}))} + \sum_{n=0}^{N-1} \textcolor{blue}{\frac{1}{2}} u_n^T R u_n  \\
        &+ \lambda_0^T \Big[ p^0 + D_1 L_d(q^0, q_1) + f_d^{-}(q^0, q_{1}, u_0) \Big] \nonumber\\
        &+ \sum_{n=1}^{N-1} \lambda_{n}^T \Big[ D_1 L_d(q_n, q_{n+1}) + D_2 L_d(q_{n-1}, q_{n}) \nonumber\\
        &+ f_{d}^{-}(q_n, q_{n+1}, u_n) + f_{d}^{+}(q_{n-1}, q_{n}, u_{n-1}) \Big]\nonumber
\end{align}
The discrete variation of the augmented objective $\delta \tilde{J}_d = 0$ has to vanish for variations $\delta u_n$, $\delta \lambda_n$ and $\delta q_n$ with boundary conditions $\delta q_0 = 0$ that is directly enforced as $q_0 = q^0$ at the initial time node is specified in problem (\ref{eqn:OCP}). The variation of the three types of variables leads to three sets of equations. The variation w.r.t the adjoint variables leads to the discrete \textrm{\textsc{Euler}}-\textrm{\textsc{Lagrange}} equations, the constraints in (\ref{eqn:OCP}). The variation with respect to the configuration variable yields the adjoint equations, reading with rearrangement of terms as follows:
\begin{subequations}\label{eqn:adj}
   \begin{alignat}{2}
        &\lambda_{N-1}^T  [ D_2 D_1 L_d(q_{N-1}, q_{N}) + D_2 f_d^-(q_{N-1}, q_N, u_{N-1})]\label{eq:adj1}\\
        &= - S_q (q_N - q^N) - \textcolor{blue}{S_p \left[ \, p_N(q_{N-1},q_N,u_{N-1}) - p^N\,\right]}\notag\\
        &\times [ D_2 D_2 L_d(q_{N-1}, q_N) + D_2 f_{d}^+(q_{N-1}, q_N, u_{N-1})]\notag\\
        \notag\\
        &\lambda_{N-2}^T  [D_2 D_1 L_d(q_{N-2}, q_{N-1}) + D_2 f_{d}^-(q_{N-2}, q_{N-1}, u_{N-2})] \label{eq:adj2}\\
        &+ \lambda_{N-1}^T \Big[  D_2 D_2 L_d(q_{N-2}, q_{N-1}) + D_2 f_{d}^+(q_{N-2}, q_{N-1}, u_{N-2})\notag\\
        &+  D_1 D_1 L_d(q_{N-1}, q_{N}) + D_1 f_{d}^-(q_{N-1}, q_{N}, u_{N-1})\Big] \notag \\
        &= \textcolor{blue}{-S_p \left[ \, p_N(q_{N-1},q_N,u_{N-1}) - p^N \right]}\notag\\
        & \textcolor{blue}{\times [D_1 D_2 L_d (q_{N-1}, q_N) + D_1 f^+_{d} (q_{N-1}, q_N, u_{N-1})]} \notag
    \end{alignat}
    \begin{alignat}{2}
        0 &= \lambda_{n-1}^T  [D_2 D_1 L_d(q_{n-1}, q_{n})+D_2 f_{d}^-(q_{n-1}, q_{n}, u_{n-1})] \label{eq:adj3}\\
        &+\lambda_{n}^T \Big[  D_2 D_2 L_d(q_{n-1}, q_{n}) + D_2 f_{d}^+(q_{n-1}, q_{n}, u_{n-1}) \notag\\
        &+ D_1 D_1 L_d(q_n, q_{n+1}) + D_1 f_{d}^-(q_{n}, q_{n+1}, u_{n}) \Big]\notag\\
        &+ \lambda_{n+1}^T  \left[ D_1 D_2 L_d(q_{n}, q_{n+1}) + D_1 f_{d}^+(q_n, q_{n+1}, u_n)\right],\quad\text{for} ~n=N-2,~...,~1 \notag
    \end{alignat}
\end{subequations}
The discrete variational principle directly provides the boundary conditions \eqref{eq:adj1} and \eqref{eq:adj2} for the two last adjoint variables, as no boundary conditions for the state variables are prescribed at these time nodes. The variation w.r.t. the input $u_n$ yields the optimality conditions. Note that the last equation is different:
\begin{subequations}\label{eqn:optU}
    \begin{align}
    0&= R u_n + \lambda^T_{n} D_3 f_d^-(q_n, q_{n+1}, u_n)\\
     &+ \lambda^T_{n+1} D_3 f_d^+(q_n, q_{n+1}, u_n),\qquad\text{for}~n=0,~...,~N-2 \nonumber\\
     &\nonumber\\
    0&= R u_{N-1} + 
    \lambda^T_{N-1} D_3 f_{d}^-(q_{N-1}, q_{N}, u_{N-1})\\
    & \textcolor{blue}{+ S_p \left[ \, p_N(q_{N-1},q_N,u_{N-1}) - p^N \right] D_3 f_{d}^+(q_{N-1}, q_{N}, u_{N-1})}\nonumber
    \end{align}
\end{subequations}
The discrete \textrm{\textsc{Euler}}-\textrm{\textsc{Lagrange}} equations (\ref{eqn:DEL}) can be solved forward in time and the adjoint equations (\ref{eqn:adj}) backward in time sequentially given the configuration path to determine the discrete adjoint variables as a shooting method while using the input equations (\ref{eqn:optU}) to update the input. \textcolor{blue}{Such a direct shooting algorithm directly uses the equations derived above and thus is simple to implement. However, an appropriately small time step $h$ is necessary for stable integration in both directions in time.} \textcolor{blue}{The discrete optimization problem with respect to $q_d$, $u_d$ and $\lambda_d$ can also be solved by applying an interior point algorithm \cite{Waechter2005} or sequential quadratic programming \cite{Gill2005}. In those, the variational integrator is used as equality constraints for the optimization as in \eqref{eqn:OCP}}

\subsection{Application of the Discrete Adjoint Method to a Mathematical Pendulum}\label{subsec:pend}
\begin{figure}[ht]
    \centering
%
%
%
%
%
%
%
%
    \includegraphics{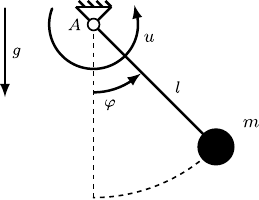}
    \caption{Torque-controlled mathematical pendulum.}
    \label{fig:pendulum}
 \end{figure}

\begin{figure}
\centering
	\begin{subfigure}{0.45\textwidth}
	  \centering
\centering\resizebox{0.9\textwidth}{!}{\includegraphics{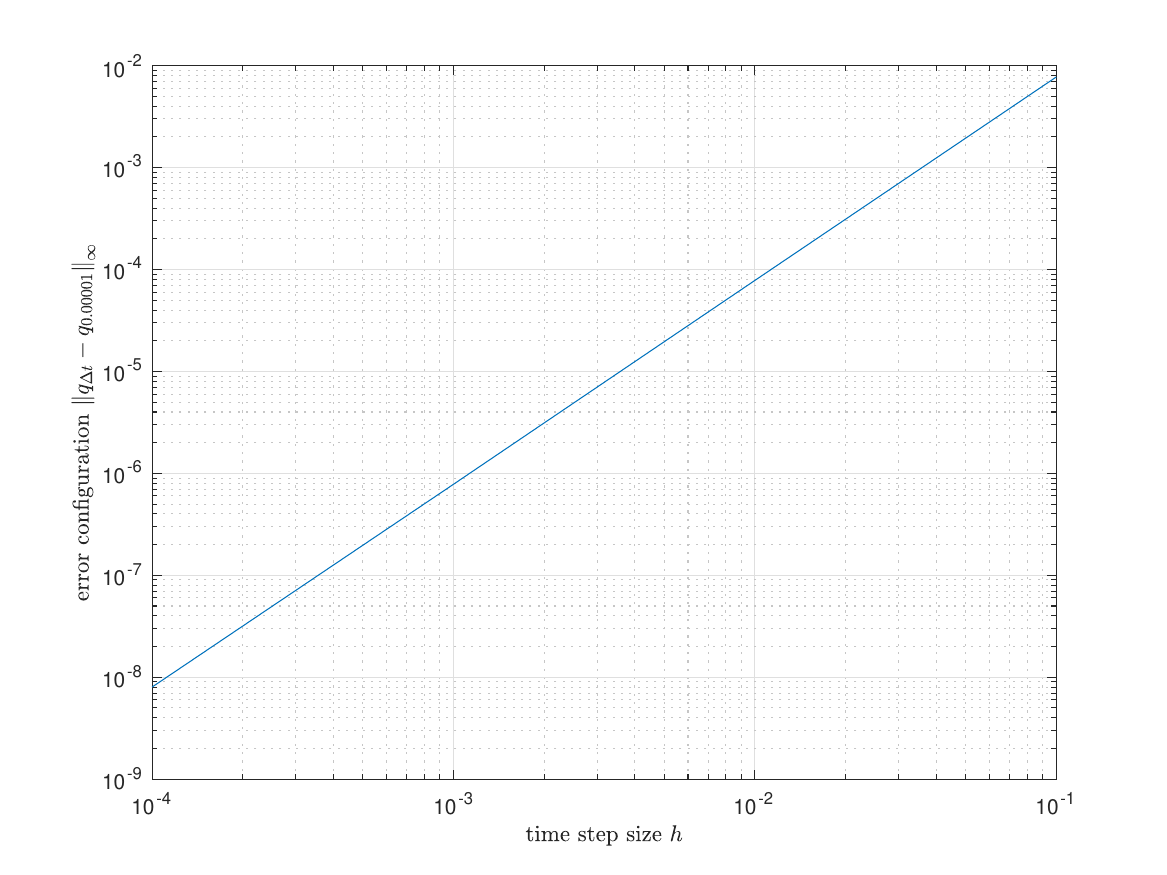}}
    \caption{Error of configuration $q_d$ versus time step for the mathematical pendulum in minimal coordinates.}\label{fig:econf}
	\end{subfigure}
 \quad
	\begin{subfigure}{0.45\textwidth}
	  \centering
	\centering\resizebox{0.9\textwidth}{!}{\includegraphics{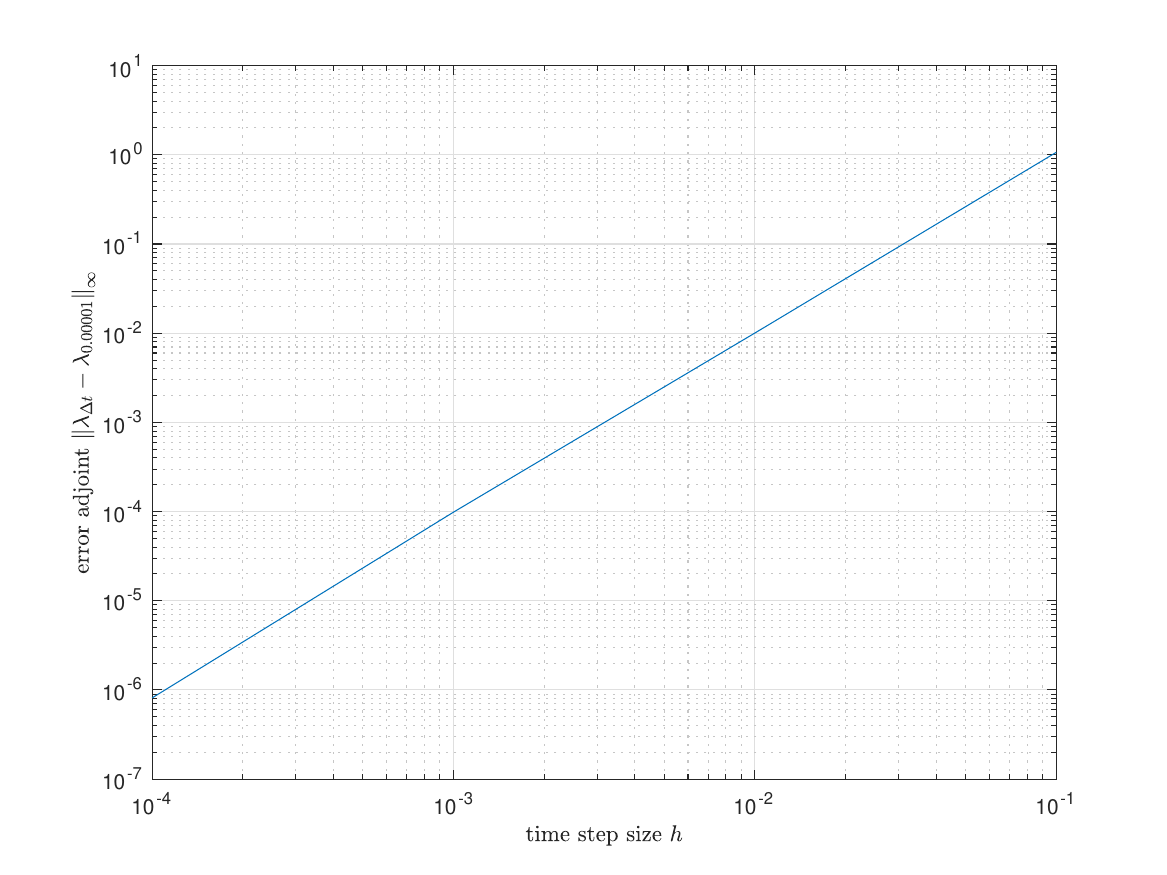}}
    \caption{Error of adjoint variables $\lambda_d$ versus time step for the mathematical pendulum in minimal coordinates.}\label{fig:eAdj}
	\end{subfigure}%
\caption{Error of the configuration $q_d$ and adjoint variable $\lambda_d$.}
\end{figure}

\begin{figure}
\centering
	\begin{subfigure}{0.5\textwidth}
	  \centering
	\resizebox{0.9\textwidth}{!}{\includegraphics{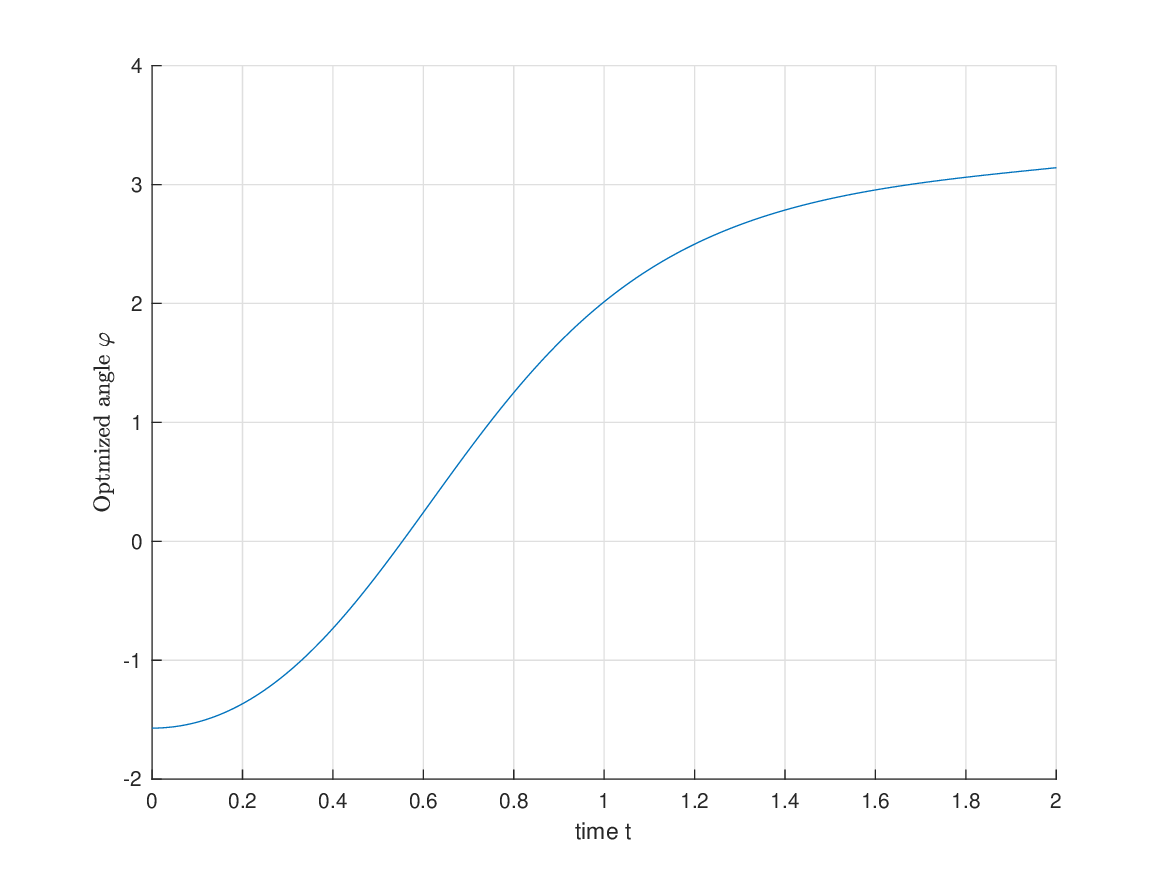}}
    \caption{Optimized angle $\varphi$.}\label{fig:ominiphi}
	\end{subfigure}%
	\begin{subfigure}{0.5\textwidth}
	  \centering
	\centering\resizebox{0.9\textwidth}{!}{\includegraphics{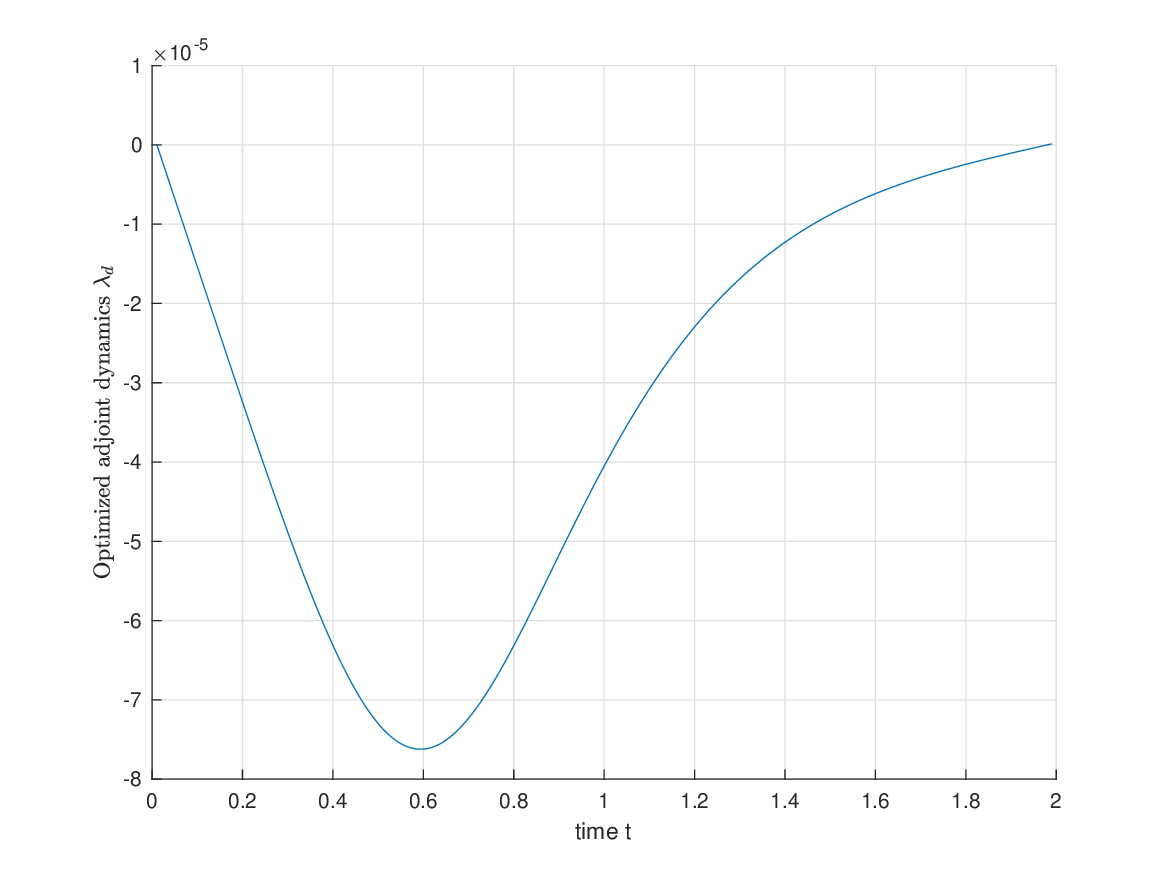}}
    \caption{Optimized adjoint variable $\lambda$.}\label{fig:ominiadj}
	\end{subfigure}
	\bigbreak
	\begin{subfigure}{0.5\textwidth}
	  \centering
  \resizebox{0.9\textwidth}{!}{\includegraphics{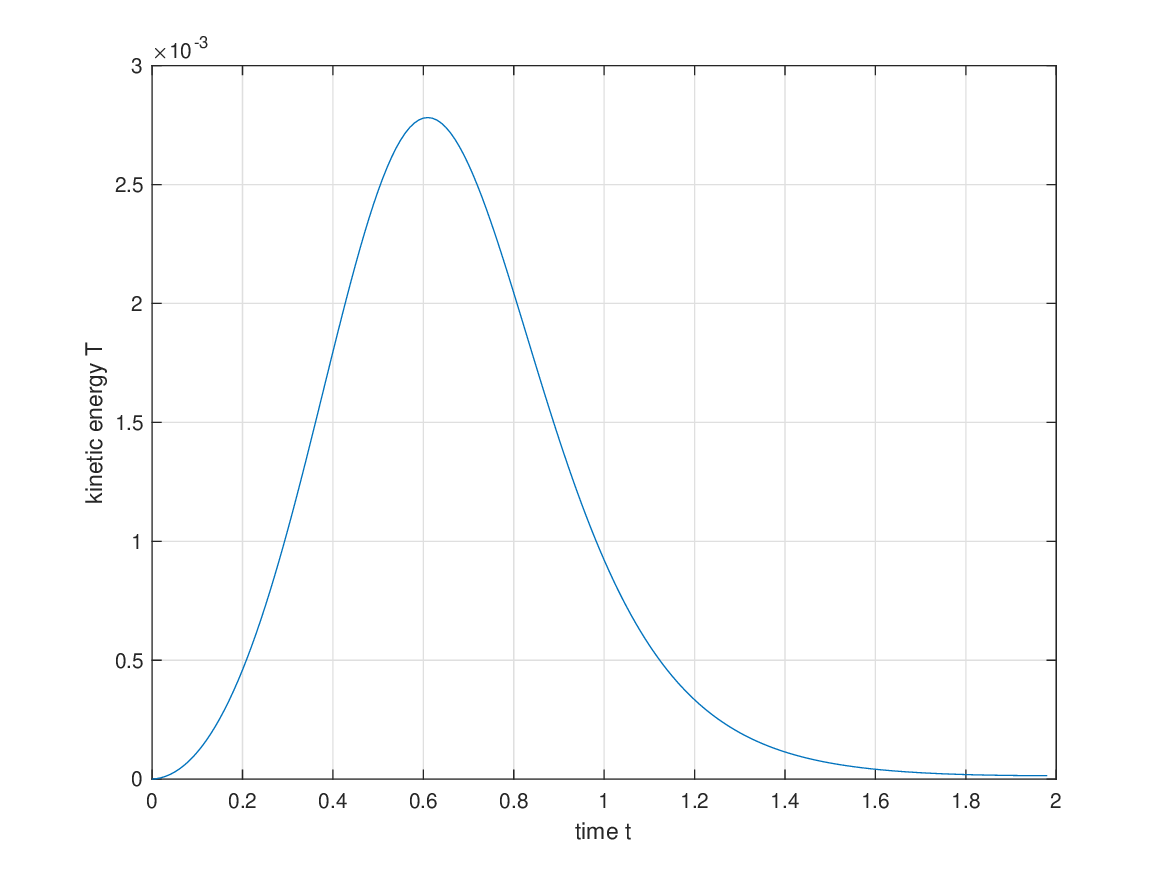}}
    \caption{Kinetic energy $T$ of the optimized upswing maneuver.}\label{fig:ominikin}
	\end{subfigure}%
	\begin{subfigure}{0.5\textwidth}
	  \centering
	    \centering\resizebox{0.9\textwidth}{!}{\includegraphics{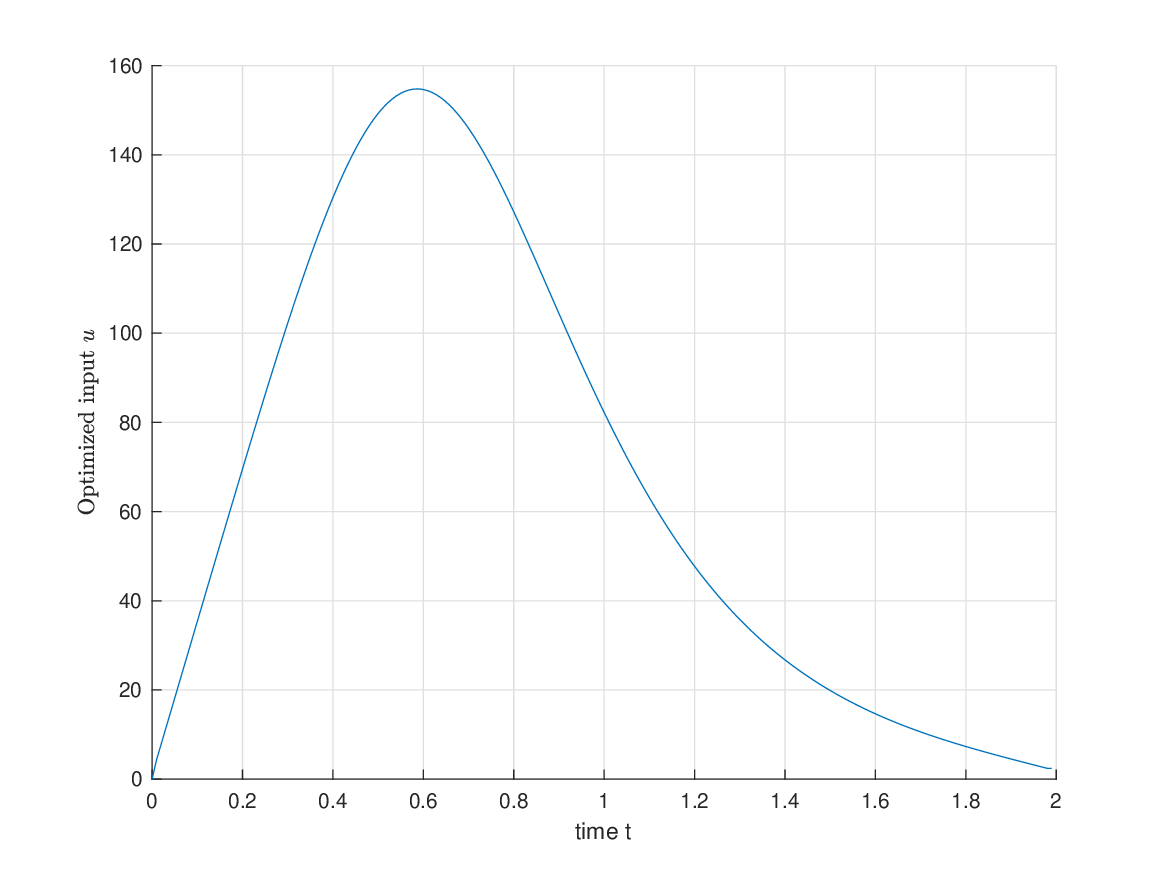}}
    \caption{Optimized input $u$.}\label{fig:ominimom}
	\end{subfigure}%
	\caption{Optimization results for the pendulum, using the discrete adjoint method and single shooting.}
	\label{fig:opt_unconst}
\end{figure}

Let us consider a mathematical pendulum as depicted in Figure \ref{fig:pendulum}, in minimal coordinates $q=\varphi$ with the Lagrangian $\mathcal{L}(\varphi, \dot{\varphi}) = \frac{1}{2} m l^2 \dot{\varphi}^2 - m g l \cos(\varphi)$ that is actuated by a torque $f=\textcolor{blue}{u}$. The discrete Lagrangian approximated with the midpoint rule is $L_d(\varphi_n, \varphi_{n+1}) = \frac{1}{2} h m l^2 (\frac{\varphi_{n+1} - \varphi_n}{h})^2 - h m g l \cos(\frac{\varphi_{n+1} + \varphi_n}{2})$ with the time step $h$.
The \textcolor{blue}{discrete} forces are $\textcolor{blue}{f_d^{\pm}(\varphi_n, \varphi_{n+1},u_n) = \frac{1}{2} h u_n}$. 
\textcolor{red}{We wish to perform an optimal upswing maneuver. Thus, the initial configuration and momentum are $\phi^0 = 0$ and $p^0 = 0$, and t}\textcolor{blue}{he desired end configuration is $q^N = \varphi^N = \pi$. The end momentum has to vanish $p^N = 0$.} The first slot derivatives of the discrete Lagrangian used for the discrete  \textrm{\textsc{Euler}}-\textrm{\textsc{Lagrange}} equations are:
\begin{align}
    D_1 L_d(\varphi_n, \varphi_{n+1}) &= - m l^2 \frac{\varphi_{n+1}-\varphi_n}{h} + \frac{h}{2} \, m g l \, \sin\left(\frac{\varphi_{n+1}+\varphi_n}{2}\right)\\
    D_2 L_d(\varphi_{n-1}, \varphi_n) &= \hphantom{-} m l^2 \frac{\varphi_{n}-\varphi_{n-1}}{h} + \frac{h}{2} \, m g l \, \sin\left(\frac{\varphi_{n}+\varphi_{n-1}}{2}\right)
\end{align}
The time stepping scheme (\ref{eqn:DEL}) for the configuration is:
\begin{equation}\label{eqn:timestep}
    \begin{aligned}
        0=\frac{\varphi_{n+1} - 2 \varphi_n + \varphi_{n-1}}{h} &- \frac{h}{2} \frac{g}{l} \, \text{sin}\left(\frac{\varphi_{n+1} + \varphi_{n}}{2}\right)\\
        & - \frac{h}{2} \frac{g}{l} \, \text{sin}\left(\frac{\varphi_n + \varphi_{n-1}}{2}\right) - h\frac{\textcolor{blue}{u}_n + \textcolor{blue}{u}_{n-1}}{2}
    \end{aligned}
\end{equation}
It is initialized with 
\begin{equation}
    0 = p^0 - m l^2 \frac{\varphi_{1}-\varphi_0}{h} + \frac{h}{2} \, m g l \, \text{sin}\left(\frac{\varphi_{1}+\varphi_0}{2}\right) + h \frac{\textcolor{blue}{u}_0}{2}
\end{equation}
The second derivatives of the \textcolor{blue}{discrete} Lagrangian inserted in the adjoint equations (\ref{eq:adj3}) leads to:
\begin{equation}\label{eqn:adjstep}
    \begin{aligned}
        0= & \frac{\lambda_{n-1}^T - 2 \lambda_{n}^T + \lambda_{n+1}^T}{h}\\
        &- \frac{\lambda_{n}^T + \lambda_{n-1}^T}{2} \, \frac{h}{2} \frac{g}{l} \, \cos\left(\frac{\varphi_n + \varphi_{n-1}}{2}\right)\\
        &- \frac{\lambda_{n+1}^T +\lambda_{n}^T}{2} \, \frac{h}{2} \frac{g}{l} \, \cos\left(\frac{\varphi_{n+1} + \varphi_{n}}{2}\right)
    \end{aligned}
\end{equation}
Two equations according to (\ref{eq:adj1}) and (\ref{eq:adj2}) are necessary to initialize the backward integration (\ref{eqn:adjstep}) in time:
\begin{subequations}\label{eqn:adjb}
    \begin{align}
        0 &= \lambda_{N-1}^T  \left[ - \frac{m l^2}{h}  + \frac{h}{4} \, m g \, \cos\left(\frac{\varphi_{N}+\varphi_{N-1}}{2}\right) \right] + S_q (\varphi_N - \pi)\\
          &+  S_p \textcolor{blue}{\left[ m l^2 \frac{\varphi_{N}-\varphi_{N-1}}{h} + \frac{h}{2} \, m g l \, \sin\left(\frac{\varphi_{N}+\varphi_{N-1}}{2}\right) \right]}\notag\\
          &\times \left[ \frac{m l^2}{h}  + \frac{h}{4} \, m g \, \cos\left(\frac{\varphi_{N}+\varphi_{N-1}}{2}\right)\right] \notag
    \end{align}
    \begin{align}
        0 &= \frac{2\lambda^T_{N-1} - \lambda^T_{N-2}}{h}\\
          &+ \frac{\lambda_{N-1}^T + \lambda_{N-2}^T}{2} \, \frac{h}{2} \, \frac{g}{l} \cos\left(\frac{\varphi_{N-1}+\varphi_{N-2}}{2}\right) + \lambda_{N-1}^T \, \frac{h}{4} \frac{g}{l} \, \cos\left(\frac{\varphi_{N}+\varphi_{N-1}}{2}\right) \notag\\
          &+  S_p \textcolor{blue}{\left[ m l^2 \frac{\varphi_{N}-\varphi_{N-1}}{h} + \frac{h}{2} \, m g l \, \sin\left(\frac{\varphi_{N}+\varphi_{N-1}}{2}\right) \right]} \notag\\
          &\textcolor{blue}{\times \left[\frac{m l^2}{h}  + \frac{h}{4} \, m g \, \cos\left(\frac{\varphi_{N}+\varphi_{N-1}}{2}\right)\right]}\notag
    \end{align}
\end{subequations}
The equations for the input are:
\begin{subequations}\label{eqn:input}
    \begin{align}
        0&=R h \, \textcolor{blue}{u}_n + h \frac{\lambda^T_{n} + \lambda^T_{n+1}}{2},\quad\text{for}~n=0,~...,~N-2 \label{eqn:input_1}\\
        0&=R h \, \textcolor{blue}{u}_{N-1} + \frac{h}{2} \lambda^T_{N-1}\label{eqn:input_2}\\
        &\textcolor{blue}{+ \, \frac{h}{2} S_p \textcolor{blue}{\left[ m l^2 \frac{\varphi_{N}-\varphi_{N-1}}{h} + \frac{h}{2} \, m g l \, \sin\left(\frac{\varphi_{N}+\varphi_{N-1}}{2}\right) \right]}}.\nonumber
    \end{align}
\end{subequations}
The convergence of the configuration $\textcolor{blue}{q_d}$ and the adjoint variables $\textcolor{blue}{\lambda_d}$ is illustrated in Figures \ref{fig:econf} and \ref{fig:eAdj}, respectively. A simulation time of $T=2$ and a constant input of $u^n=1$, for $n=0,~...,~N-1$ is used; the pendulum has a length of $L=1$ with a gravitational constant of $g=9.81$. \textcolor{red}{The mass of the pendulum is $m=1$. For the input weight $R=10^{-5} h$ is used. The weights in the \textrm{\textsc{Mayer}} term are $S_q = 10^3$ and $S_p = 10^{-2}$. These values were chosen to obtain solutions that achieve the upswing of the pendulum to the upper equilibrium point, with minimal effort. Larger values for the input weighting lead to solutions with end configuration at the lower equilibrium of the pendulum. }The absolute error in these plots is computed using the infinity norm of the difference of the variables and a reference solution, \textcolor{blue}{$(q_{\mathrm{ref}}, \lambda_{\mathrm{ref}})$}, which is a simulation with a fine discretization of \textcolor{blue}{$h=10^{-5}$}, \textcolor{blue}{$\left\Vert q_d - q_{\mathrm{ref}} \right\Vert_{\infty}$} and $\textcolor{blue}{\left\Vert \lambda_d - \lambda_{\mathrm{ref}} \right\Vert_{\infty}}$, respectively. The convergence rate for the configuration and adjoint variables is equal, we observe second order convergence. This is in accordance to the theoretical results in \cite{Marsden2011}. 
These convergence results are derived for the forward integration of the time stepping scheme (\ref{eqn:timestep}) and the subsequent backwards solution of (\ref{eqn:adjstep}) using the configuration variables calculated with the same time step width.\\
\\
\textcolor{blue}{The optimized motion of the pendulum is depicted in Figures \ref{fig:ominiphi}, \ref{fig:ominiadj}, \ref{fig:ominikin} and \ref{fig:ominimom}. The momentum $p$ and the kinetic energy $T$ is close to zero at the end of the simulation with the optimized input acting on the pendulum.}\\
\\
\\
REMARK 2: \textrm{\textsc{Pontryagin}}'s maximum \textcolor{red}{principle} leads to \textcolor{blue}{necessary conditions for optimality in the continuous-time setting. The resulting adjoint equations are $\ddot{\lambda}^T - g/l \lambda \cos \varphi = 0$ and the control equations are $R u + \lambda = 0$. It can be checked that the discrete equations \eqref{eqn:adjstep} and \eqref{eqn:input_1} are the corresponding discrete versions of these equations when discretized using a midpoint rule. The discrete boundary conditions \eqref{eqn:adjb} and \eqref{eqn:input_2}, however, are not so easy to relate to their continuous counterparts. We plan to address this very issue in a future publication.}

\section{Discrete Adjoint Method for Variational Integration of Constrained Dynamics}

\subsection{Variational Integration of Constrained Dynamics}\label{sec:VIC}%
The derivation of variational integrators for constrained systems that use null space projection and nodal reparametrization \cite{Leyendecker2008} is shortly summarized in the following section, using similar steps as in Section \textcolor{blue}{\ref{sect:DAMVIs}}. The discrete adjoint method for such systems is derived thereafter \textcolor{blue}{similar to Section \ref{sec:DAMVI}}.\\
Up \textcolor{blue}{until} now, \textcolor{blue}{we have workd in local coordinates directly on the configuration manifold $\mathcal{Q}$}. However, it \textcolor{blue}{can be} advantageous to consider \textcolor{blue}{$\mathcal{Q}$ an ambient (vector) space parametrized by} redundant coordinates, and \textcolor{blue}{constrain the motion by constraints. Given a scleronomic, holonomic constraint function $g: \mathcal{Q} \to \mathbb{R}^m$, the constraint submanifold is then}
\begin{equation}
    \mathcal{M} := \{q \in \mathcal{Q} ~|~ g(q)=0\}.
\end{equation}
We assume that the Jacobian $\frac{\partial g}{\partial q}$ has full rank $m$, so the dimension of the constraint manifold is $n-m$, the number of degrees of freedom of the mechanical system. \textcolor{blue}{We also assume consistent initial conditions $(q^0,\dot{q}^0)$ that fulfill the constraints on configuration level $g(q^0)=0$ as well as on velocity level $\frac{d}{dt} g(q^0) = \frac{\partial g(q^0)}{\partial q} \dot{q}^0=0$.}

A \textrm{\textsc{Lagrange}} multiplier $\nu$ \textcolor{blue}{is used to enforce the} constraint \textcolor{blue}{by appending the term $- g(q)^T \nu$} to the Lagrangian in the action integral.  \textcolor{blue}{Thus, t}he \textrm{\textsc{Lagrange}}-\textrm{\textsc{d'Alembert}} principle in this setting reads
\begin{equation}
    0 = \delta \int_0^T \left[ \mathcal{L} \big(q(t), \dot{q}(t)\big)  - g(q)^T \nu \right] ~dt + \int_0^T f_\mathcal{L} \big(q(t), \dot{q}(t), u(t) \big) \delta q~ dt,\quad \forall \delta q, \delta \nu
\end{equation}
with $\delta q(0) = \delta q(T) = 0$. The constraint part of the action integral is approximated with the trapezoidal rule:
\begin{equation}
    \int_{t_n}^{t_{n+1}} \textcolor{blue}{g(q(t))^T \nu(t)} ~dt \approx \frac{1}{2} \big[  g_d(q_{n}) \nu_n + g_d(q_{n+1}) \nu_{n+1} \big]
\end{equation}
with $g_d(q_n) = h g(q_n)$. Including this in the discrete variational principle in (\ref{eqn:LAprinc}), in the constrained case, the discrete variational principle the variation of the discrete action sum with the variations $\delta q_n$ and $\delta \nu_n$ and $\delta q_0 = \delta q_N = 0$ with subsequent rearrangement of terms leads to the discrete, constrained \textrm{\textsc{Euler}}-\textrm{\textsc{Lagrange}} equations
\begin{subequations}
	\begin{align}
	        0 &= D_1 L_d(q_n, q_{n+1}) + D_2 L_d(q_{n-1}, q_n) + \frac{\partial g_d(q_n)}{\partial q_n}^T \nu_n\label{eq:FDEL_with_constraints}\\
	          &+ f_d^-(q_n, q_{n+1}, u_n) + f_d^-(q_{n-1}, q_n, u_{n-1}) \nonumber\\
	        0 &= g(q_{n+1})
	\end{align}
\end{subequations}
of dimension $n+m$. To reduce the dimension \textcolor{blue}{of \eqref{eq:FDEL_with_constraints} from $n$} to $n-m$ and \textcolor{blue}{eliminate the \textrm{\textsc{Lagrange}} multipliers, thus avoiding} conditioning problems \textcolor{blue}{related to these}, a discrete null space matrix \textcolor{blue}{$P(q_n) \in \mathbb{R}^{n \times (n-m)}$, with columns spanning the tangent space $T_{q_n}\mathcal{M}$,} can be applied that only depends on quantities at the current step such that the constraint forces are eliminated. Further, a nodal reparametrization $q_{n+1} = F_d(q_n, v_{n+1})$ with $v_{n+1} \in \mathcal{V} \subseteq \mathbb{R}^{n-m}$ is then used to eliminate the constraints as $g(F_d(q_n, v_{n+1})) = 0$, $\forall v_{n+1} \in \mathcal{V},~\text{for}~n=0,~...,~N-1$. \textcolor{blue}{Together with the null space matrix,} the reparametrization $F_d : \mathcal{V} \times \mathcal{Q} \to \mathcal{M}$ leads to the integration scheme 
\begin{align}\label{eqn:DElc}
        P^T(q_n) &\big[ D_1 L_d\big(q_n, F_d(q_n, v_{n+1}) \big)  + D_2 L_d\big(q_{n-1}, q_n \big) \\	
        &+ f_d^-\big(q_n,  F_d(q_n, v_{n+1}), u_n \big)\nonumber\\
        &+ f_d^+ \big(q_{n-1}, q_n, u_{n-1} \big) \big] =0,\quad\text{for}~n=1,~...,~N-1\nonumber
\end{align}
that has to be iteratively solved for $v_{n+1}$ in each time step, given $q_{n-1}, q_n, u_{n-1}$ and $u_n$.\\
The redundant control forces $f\textcolor{blue}{(q,u)} = B^T(q) \tau(u) \in \textcolor{blue}{\mathbb{R}^n}$ depend on the generalized control forces \textcolor{blue}{$\tau(u) \in \mathbb{R}^{n - m}$} and the input transformation matrix $B^T(q) \in \mathbb{R}^{n \times (n-m)}$ that must be chosen such that the consistency with the constraints and consistency of momentum maps is ensured \cite{Leyendecker2011}. The discrete approximations of the \textcolor{blue}{redundant} forces
\textcolor{blue}{\begin{align*}
f_d^-(q_n,q_{n+1},u_n) &= \frac{h}{2} B^T(q_n) \tau(u_n)\\
f_d^+(q_n,q_{n+1},u_n) &= \frac{h}{2} B^T(q_{n+1}) \tau(u_n)
\end{align*}}
capture the effect of the generalized forces acting on the time $[t_n, t_{n+1}]$. \textcolor{blue}{We have assumed that} $u$ is approximated constant in each time interval.

\subsection{Derivation of the Discrete Adjoint Method for Variational Integration of Constrained Dynamics}\label{sec:DAMVIC}%
The constrained setting with null space projection and reparametrization for a mechanical system leads to implicit equations of minimal dimension. The discrete adjoint method applied to such a system leads to adjoint variables of minimal dimension $n-m$. It also involves the null space projection for the adjoint equations.\\
The starting point is a problem such as in equation (\ref{eqn:OCP}), but now constrained by the discrete \textrm{\textsc{Euler}}-\textrm{\textsc{Lagrange}} equations for the constrained system with null space projection and nodal reparametrization (\ref{eqn:DElc}) as in \cite{Leyendecker2011}. Similar to the procedure outlined in Section \ref{sec:DAMVI}, the objective is augmented with the discrete \textrm{\textsc{Euler}}-\textrm{\textsc{Lagrange}} equations. As these equations are defined on $\mathcal{M}$ using the nodal reparametrization, $q_{n+1} = F_d(q_n, v_{n+1})$, \textcolor{red}{the adjoint variables are of the same dimension as $v_{n+1}$}.\\
An objective $J_d$ consisting of a \textrm{\textsc{Mayer}} term and an integral term quadratic in the control, similar \textcolor{blue}{to} the discrete adjoint method for systems without constraints in Equation (\ref{eqn:OCP}) is considered\textcolor{blue}{:
\begin{equation}
    J_d = \frac{1}{2} (q_N - q^N)^T S_q (q_N - q^N) + \sum_{n=0}^{N-1} \frac{1}{2} u_n^T R u_n.
\end{equation}}
However, to simplify \textcolor{blue}{matters}, the \textrm{\textsc{Mayer}} term of the momentum \textcolor{blue}{has been} omitted \textcolor{blue}{since} it can be handled \textcolor{blue}{similarly as in the unconstrained case}. The variation of the objective\textcolor{red}{, $\delta J_d$,} with respect to all variables $\delta \lambda_n$, $\delta u_n$ and $\delta v_{n+1}$ at all time steps has to vanish. The \textcolor{blue}{variation of} the redundant configuration~$\delta q_n$ \textcolor{blue}{with respect to the minimal} coordinate~$\delta v_n$ \textcolor{blue}{reads
\begin{equation}
    \delta q_n = D_2 F_d(q_{n-1}, v_n) ~\delta v_n.
\end{equation}}

The Jacobian matrix $\frac{\partial F_d}{\partial v_n}$ is a null space matrix \cite{Betsch2006}. \textcolor{blue}{After applying this relation, the adjoint equations become}
\begin{subequations}
	\begin{align}
	        &\lambda_{N-1}^T P^T(q_{N-1}) \left[ D_2 D_1 L_d(q_{N-1}, F_d(q_{N-1}, v_{N})) \right]\\
	        &= - S_q (q_N - q^N) ~D_2 F_d(q_{N-1}, v_N)\nonumber\\
	        \nonumber\\
	        &\lambda_{N-2}^T P^T(q_{N-2}) \Big[ D_2 D_1 L_d(q_{N-2}, F_d(q_{N-2}, v_{N-1})) \Big]\\
	        &= - \left\lbrace \lambda_{N-1}^T D_1 P^T(q_{N-1}) \Big[ D_1 L_d(q_{N-1}, F_d(q_{N-1}, v_{N}))+D_2 L_d(q_{N-2}, q_{N-1})\right.\nonumber\\
	        &+ f_{d}^{-}\textcolor{blue}{(q_{N-1}, F_d(q_{N-1}, v_{N}), u_{N-1}) + f_{d}^{+}(q_{N-2}, q_{N-1}, u_{N-2})}   \Big] \nonumber\\
	        &+ \left. \lambda_{N-1}^T P^T(q_{N-1}) \Big[ D_1 D_1 L_d(q_{N-1}, F_d(q_{N-1}, v_{N})) + D_2 D_2 L_d(q_{N-2}, q_{N-1}) \Big] \right\rbrace\nonumber\\
	        &\times D_2 F_d(q_{N-2}, v_{N-1})\nonumber
	\end{align}%
	\begin{align}
	        &\lambda_{n-1}^T P^T(q_{n-1}) \Big[ D_2 D_1 L_d(q_{n-1}, F_d(q_{n-1}, v_{n})) \Big] \\
	        &= - \left\lbrace \lambda_{n}^T D_1 P^T(q_n) \Big[ D_1 L_d(q_n, F_d(q_n, v_{n+1})) + D_2 L_d(q_{n-1}, q_{n}) \right.\nonumber\\
	        &+ \textcolor{blue}{f_{d}^{-}(q_n, F_d(q_n, v_{n+1}), u_n) + f_{d}^{+}(q_{n-1}, q_n, u_{n-1}) }  \Big]\nonumber\\
	        &+ \lambda_{n}^T P^T(q_n) \Big[ D_1 D_1 L_d(q_n, F_d(q_n, v_{n+1}))+ D_2 D_2 L_d(q_{n-1}, q_{n}) \Big]\nonumber\\
	        &+\left.\lambda_{n+1}^T P^T(q_{n+1}) D_1 D_2 L_d(q_{n}, q_{n+1}) \vphantom{\sum} \right\rbrace D_2 F_d(q_{n-1}, v_n)\,,~\text{for}~n=N-2,~...,~1.\nonumber
	\end{align}%
\end{subequations}
The variations with respect to the input variables vanish, if
\begin{subequations}
    \begin{align}
    0 &= R u_n + \lambda_{n}^T P^T(q_n)  D_3 f_{d}^{-}\textcolor{blue}{(q_n, F_d(q_n, v_{n+1}), u_n)}\\
      &+ \lambda_{n+1}^T P^T(q_{n+1})D_3 f_{d}^{+}\textcolor{blue}{(q_n, F_d(q_n, v_{n+1}), u_{n})},~\text{for}~n=1,~...,~N-2\nonumber\\
      &\nonumber\\
    0 &= R u_{N-1} + \lambda_{N-1}^T P^T(q_{N-1})  \textcolor{blue}{D_3 f_{d}^{-}(q_{N-1}, F_d(q_{N-1}, v_{N}), u_{N-1})}
    \end{align}
\end{subequations}
holds. The evaluation of these equations can be used to update the input variables in a shooting method.\\
\subsection{Discrete Adjoint Method for a Mathematical Pendulum described as Constrained System}%
\begin{figure}
\centering
	\begin{subfigure}{0.45\textwidth}
	  \centering
   \centering\resizebox{0.9\textwidth}{!}{\includegraphics{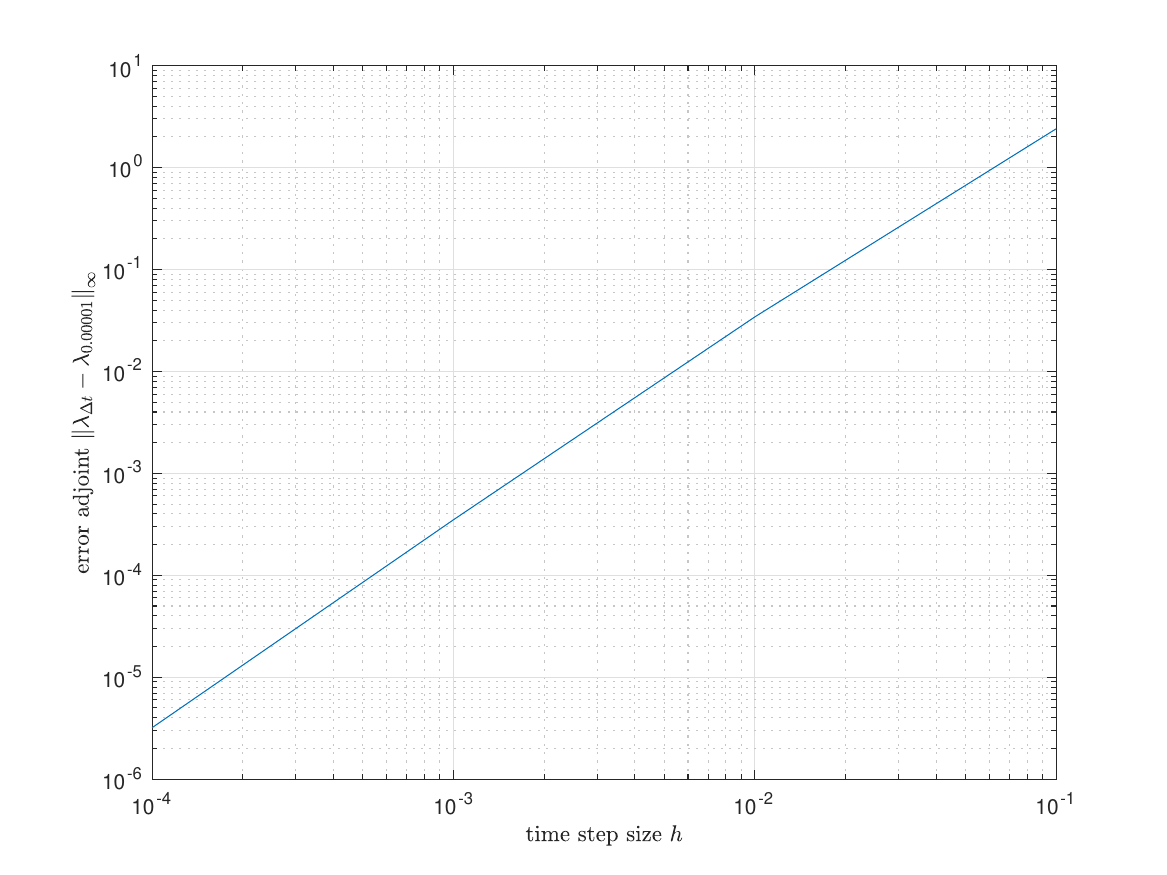}}
    \caption{Error of adjoint variables $\lambda_d$ versus time step for the mathematical pendulum as constrained system.}\label{fig:eAdjC}
	\end{subfigure}
 \quad
	\begin{subfigure}{0.45\textwidth}
	  \centering
 \centering\resizebox{0.9\textwidth}{!}{\includegraphics{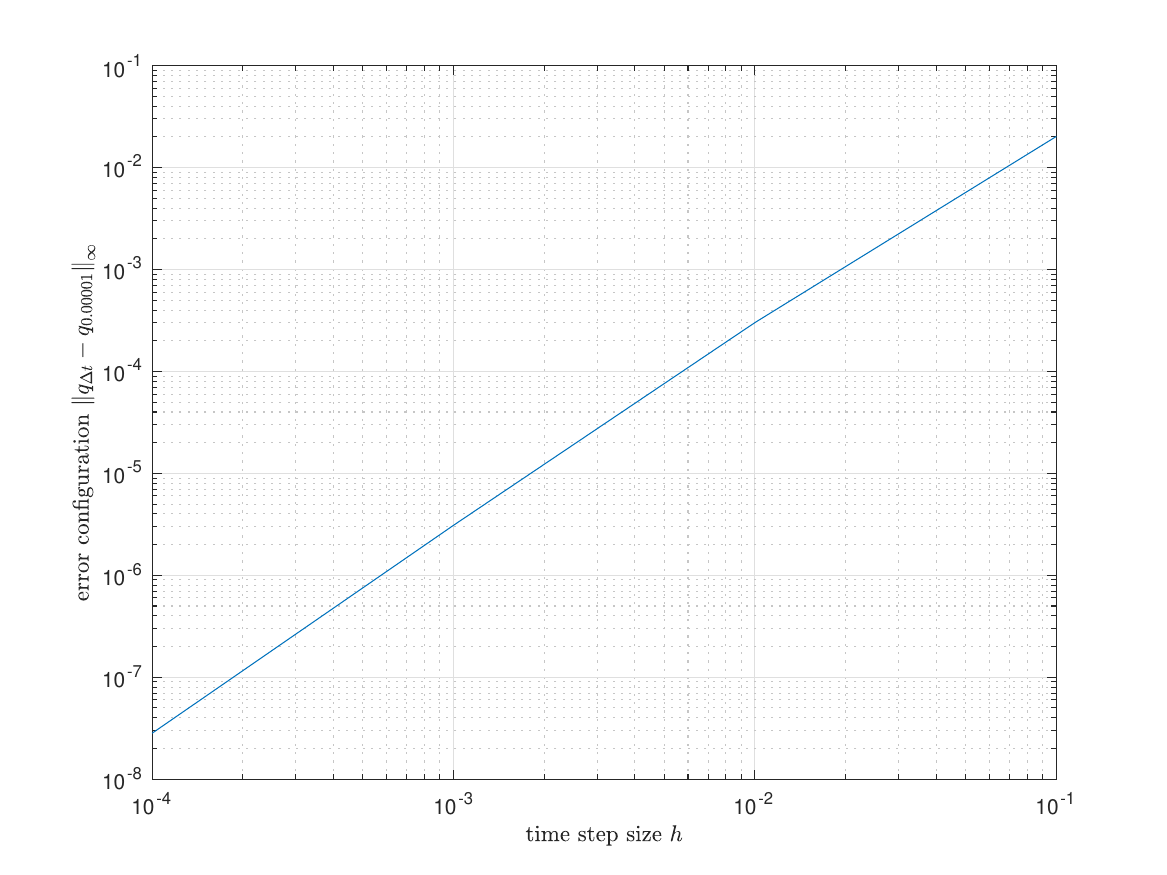}}
    \caption{Error of configuration $q_d$ versus time step for the mathematical pendulum as constrained system.}\label{fig:econfC}
	\end{subfigure}%
\caption{Error of the configuration $q_d$ and adjoint variable $\lambda_d$.}
\end{figure}

\begin{figure}
\centering
	\begin{subfigure}{0.5\textwidth}
	  \centering
	\centering\resizebox{0.9\textwidth}{!}{\includegraphics{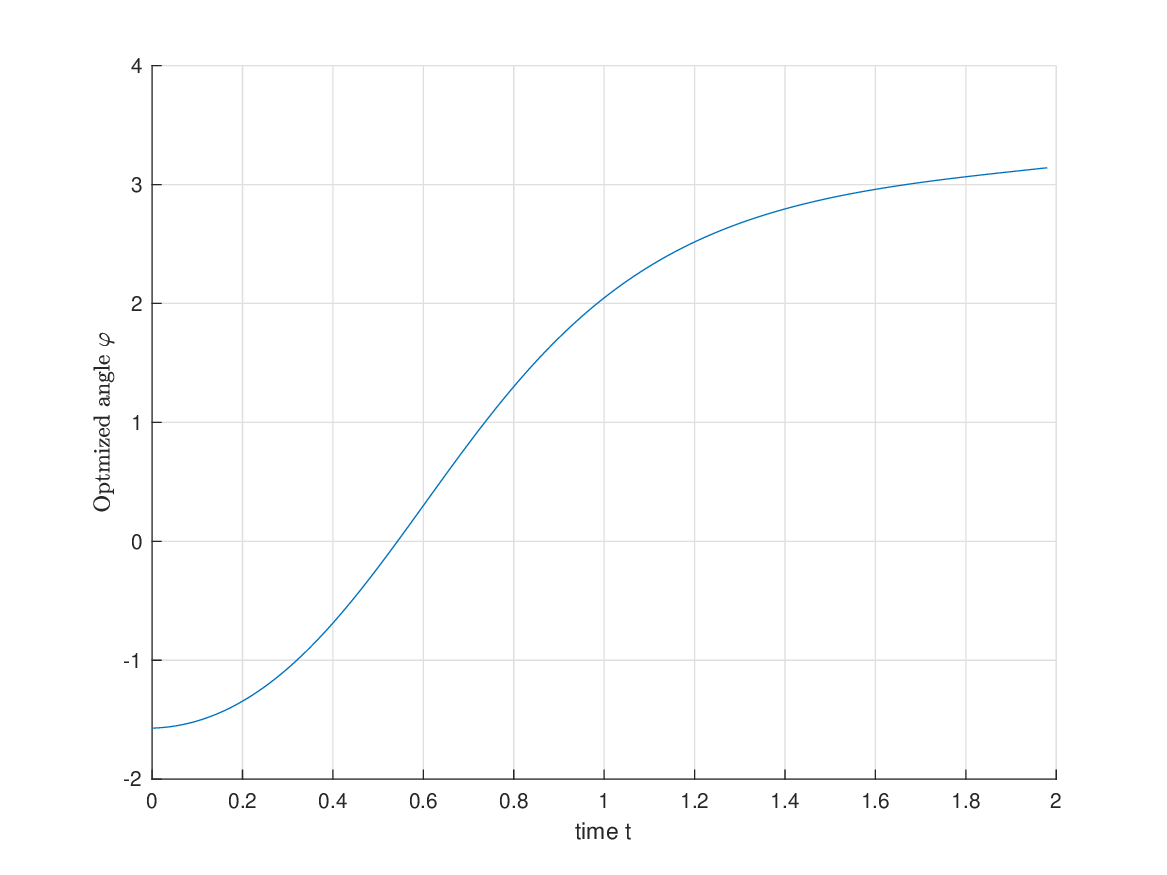}}
    \caption{Optimized angle $\varphi$, calculated as sum of $v$.}\label{fig:orephi}
	\end{subfigure}%
	\begin{subfigure}{0.5\textwidth}
	  \centering
	\centering\resizebox{0.9\textwidth}{!}{\includegraphics{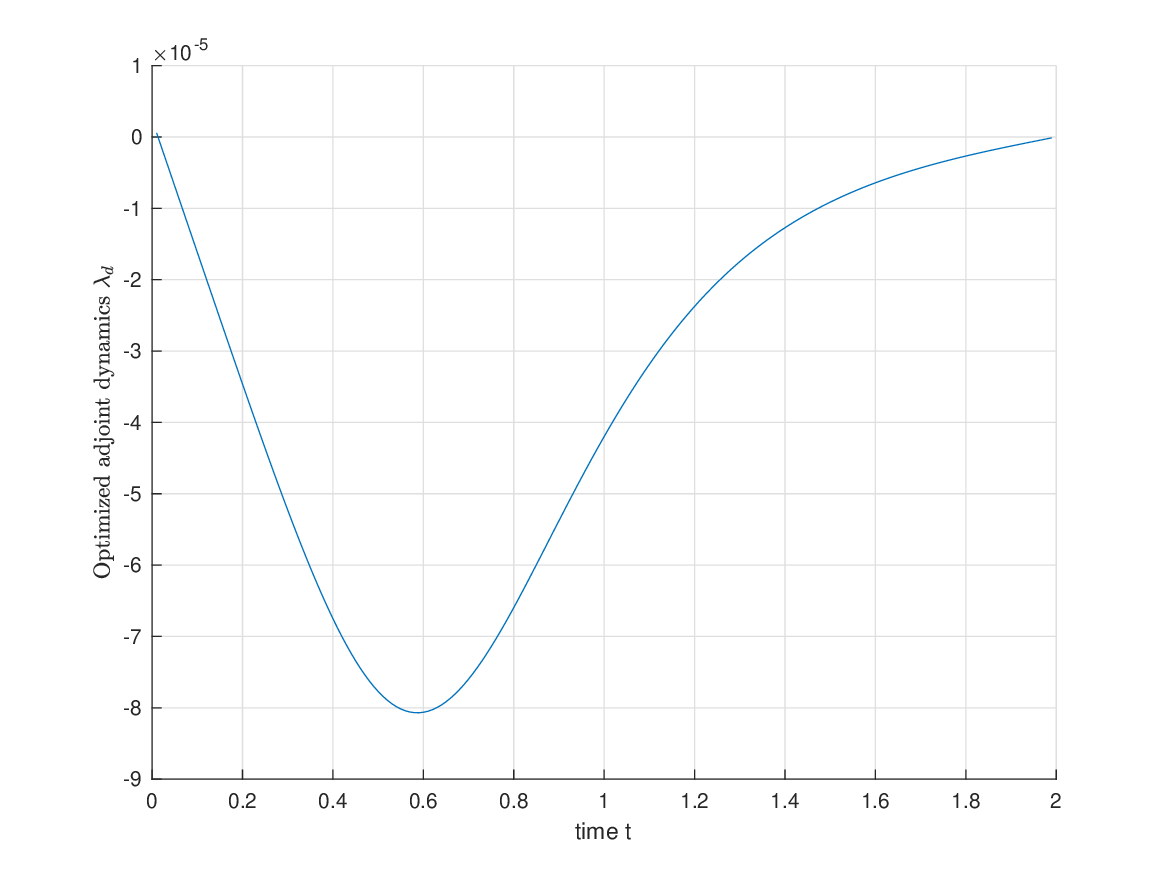}}
    \caption{Optimized adjoint variable $\lambda$.}\label{fig:oreadj}
	\end{subfigure}
	\bigbreak
	\begin{subfigure}{0.5\textwidth}
	  \centering
  \centering\resizebox{0.9\textwidth}{!}{\includegraphics{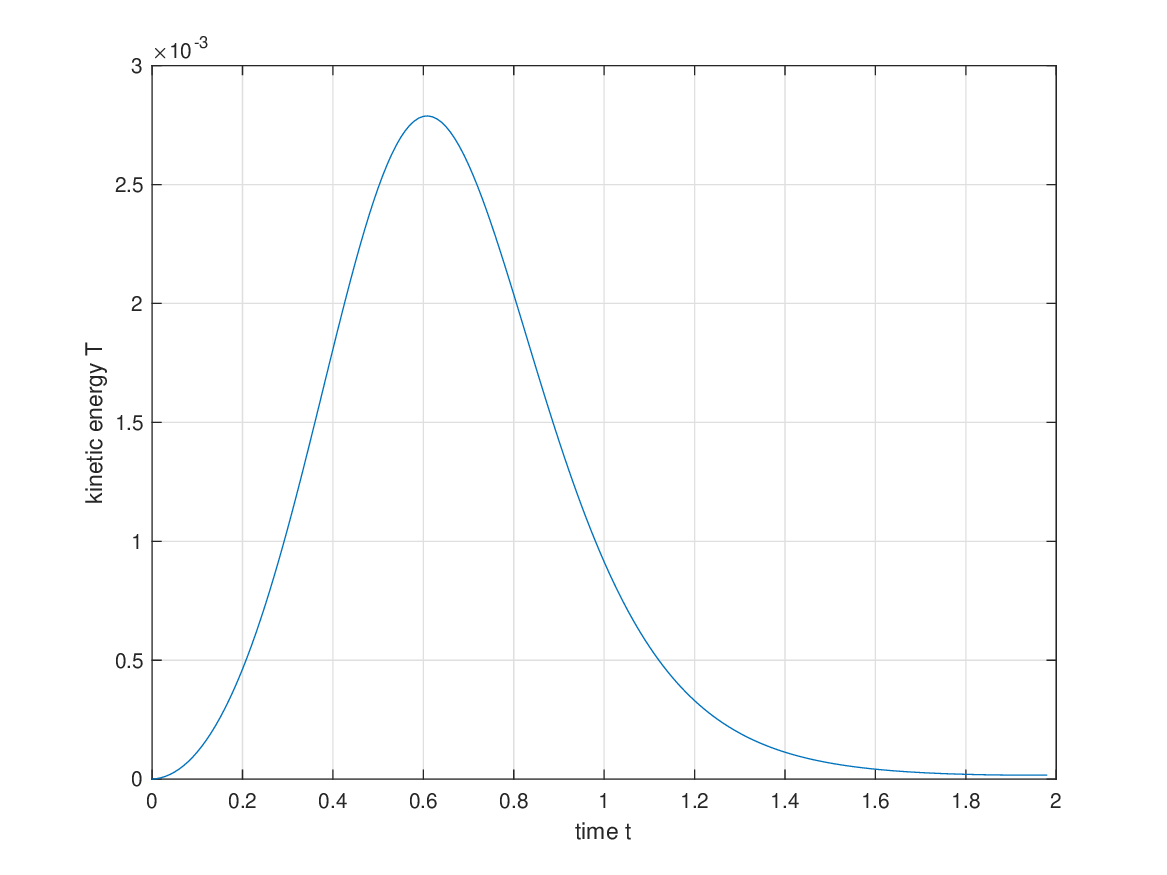}}
    \caption{Kinetic energy $T$ of the optimized upswing maneuver.}\label{fig:orekin}
	\end{subfigure}%
	\begin{subfigure}{0.5\textwidth}
	  \centering
	    \centering\resizebox{0.9\textwidth}{!}{\includegraphics{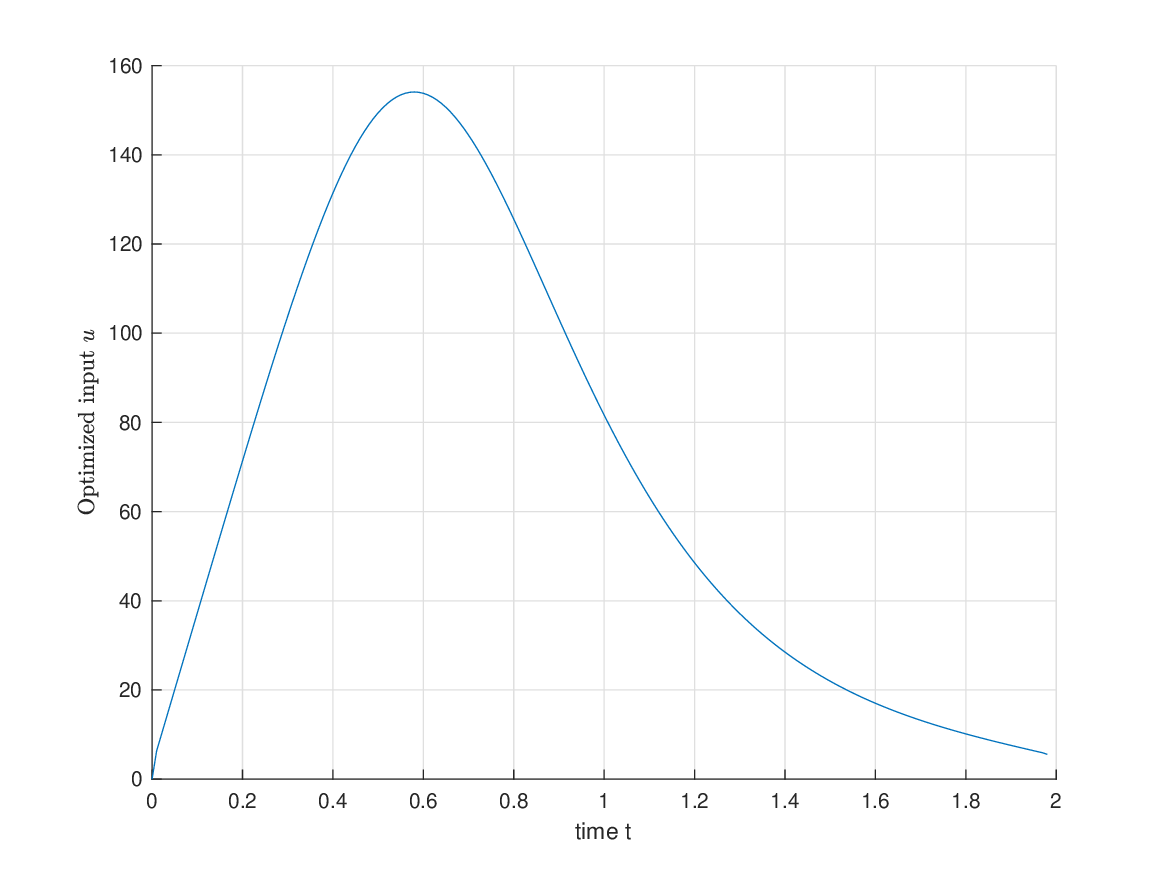}}
    \caption{Optimized input $u$.}\label{fig:oremom}
	\end{subfigure}%
	\caption{Optimization results for the pendulum in the constraint setting, using the discrete adjoint method and single shooting.}
	\label{fig:opt_constr}
\end{figure}

The mathematical pendulum is described as a constrained system \textcolor{blue}{in the ambient space $\mathcal{Q} = \mathbb{R}^2$ with redundant coordinates} $q=[x \quad y]^T$ and the constraint equation $g(q)=\textcolor{blue}{1/2(}x^2 + y^2 - l^2 \textcolor{blue}{)}$. The null space matrix is $P(q_n)^T = \textcolor{blue}{[-y_n \quad x_n]}$, the input transformation matrix is $B(q_n)^T = \textcolor{blue}{[\frac{-y_n}{2l^2} \quad \frac{x_n}{2l^2}]}$, \textcolor{blue}{the generalized force is $\tau(u) = u$,} and the nodal reparametrization reads
\begin{equation}
    q_{n+1} = F_d(q_n, v_{n+1}) = \begin{bmatrix}
    \cos(v_{n+1}) & -\sin(v_{n+1})\\
    \sin(v_{n+1}) & \cos(v_{n+1})\\
\end{bmatrix}q_n.
\end{equation}
The input variable can be interpreted as the physical torque and the variable $v$ as the incremental angle.
The Figures \ref{fig:eAdjC} and \ref{fig:econfC} show the convergence results for the pendulum in the constrained case. The adjoint variables are of minimum dimension $(n-m)$ just as the configuration variables. The error is calculated in the same way as for the unconstrained case in Section \ref{subsec:pend} as infinity norm of the difference to the reference trajectory using the same parameters. These errors are determined with solutions obtain via forward timestepping for the configuration and backward timestepping for the adjoint variables with fixed input. It can be observed in the figures that also in the constrained case the convergence rate is quadratic. However, note that the \textcolor{blue}{theoretical} results in \cite{Marsden2011} only consider the case in minimal coordinates and not the constrained case.\\
\\
\textcolor{blue}{
The optimized motion of the pendulum is depicted in Figures \ref{fig:orephi}, \ref{fig:oreadj}, \ref{fig:orekin} and \ref{fig:oremom}. The input $u$ and the kinetic energy $T$ are close to zero at the end of the simulation with the optimized input acting on the pendulum. The end configuration is weighted with $S_q=10^3$, the end momentum weight is $S_p=10^{-2}$. The weight for the input is $R=10^{-5} h$. This low weight for the input is chosen to reach the upper equilibrium position of the pendulum. It reduces the input from a constant initial guess of $1$ as well as regularizing the optimization problem.\\
The results are similar to those obtained previously by the pendulum in minimal coordinates. Small differences in the solution are visible, but show a similar optimized result.}\\

\section{Discrete Adjoint Method for Geometrically Exact Beam Dynamics}%
\textcolor{blue}{In this section, the discrete adjoint method is applied to an optimal control problem involving dynamics of a geometrically exact beam being approximated via the multisymplectic integrator found in \cite{Leitz2021}}.\\

\subsection{Geometrically Exact Beam Model}
\textcolor{blue}{The geometrically exact beam \cite{Simo1985} models a rod-like deformable body as a curve $x(t,s) \in \mathbb{R}^3$, with a rigid cross section attached to each of its points. Here $t \in [0,T]$ is used again to parametrize time, while $s \in [0, \ell]$ parametrizes the longitudinal position along the curve. The orientation of the cross section at $s$ is described by a rotation $R(t,s) \in SO(3)$. When considered as a collection of columns, $R(t,s) = [d_1(t,s), d_2(t,s), d_3(t,s)]$, the triad of vectors are known as the directors of the cross section.}\\

\begin{figure}%
    \centering%
    \includegraphics[scale=0.2]{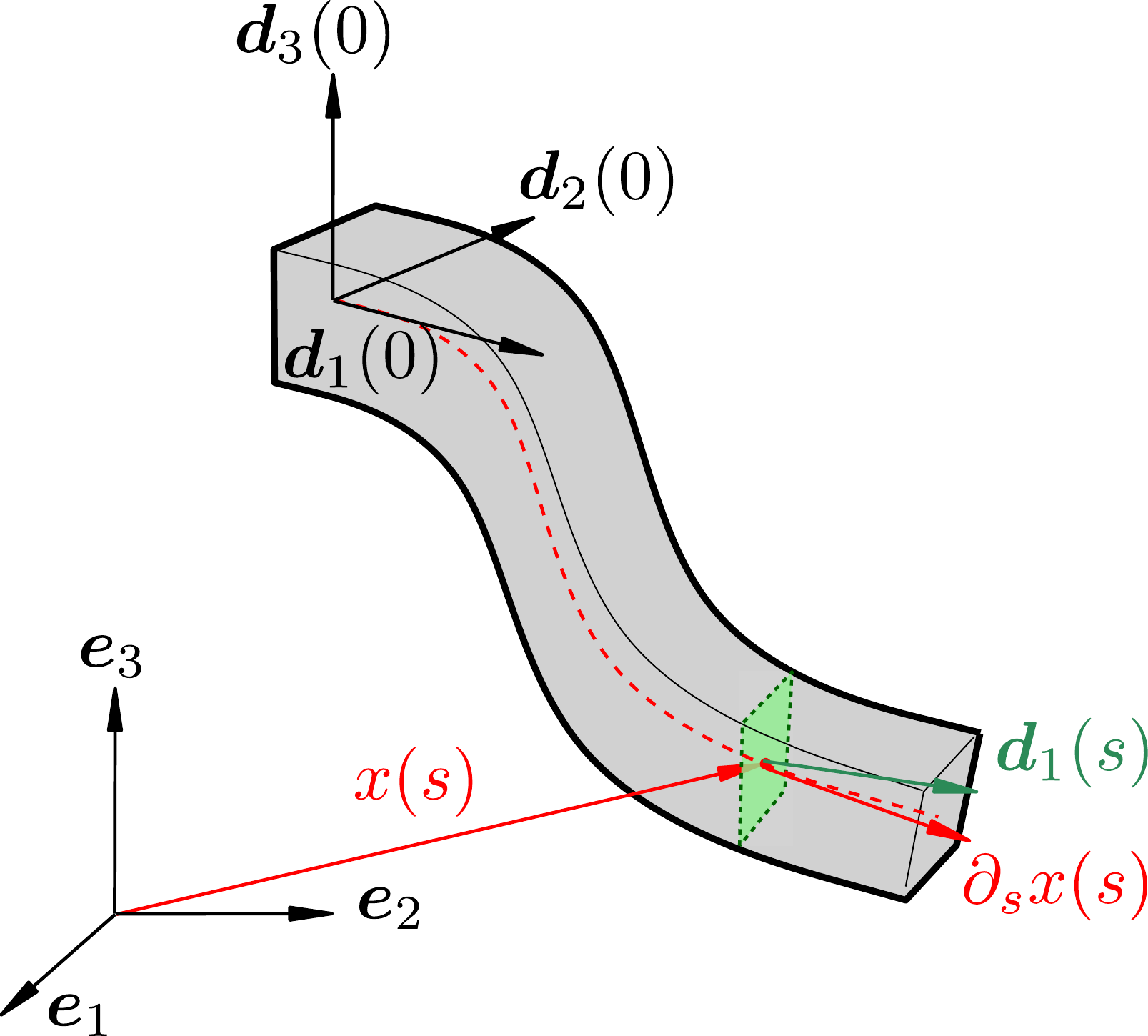}%
    \caption{Configuration of a geometrically exact beam.}
    \label{fig:beamb}
\end{figure}%

\textcolor{blue}{This can be considered as a Lagrangian field theory with configuration space $\mathcal{Q} = \mathbb{R}^3 \times SO(3)$. This space is diffeomorphic to the group of special Euclidean transformations in 3D, $SE(3)$, to which it differs only in terms of group structure. In \cite{Leitz2021,Sonneville2017}, the authors claim it to be numerically more advantageous to consider this latter space.}\\

\textcolor{blue}{If $g(t,s) = (R(t,s), x(t,s)) \in SE(3)$ denotes the configuration of a cross section, its derivatives with respect to $t$ and $s$ are related to velocities and strains respectively. More specifically,
\begin{equation*}
\begin{alignedat}{3}
(\Omega, V) &= g^{-1} \dot{g} = (R^{-1} \dot{R}, R^{-1} \dot{u}) & \; \text{body angular and linear \textcolor{red}{time derivatives}}\\
(K, W) &= g^{-1} g^\prime = (R^{-1} R^{\prime}, R^{-1} x^{\prime}) & \; \text{body angular and linear \textcolor{red}{space derivatives}}
\end{alignedat}
\end{equation*}
where we have used $\dot{X} = \frac{\partial X}{\partial t}$ and ${X}^{\prime} = \frac{\partial X}{\partial s}$ and ``body'' is meant to signify ``in the reference frame of the section itself''. Considering a reference configuration $g_{\mathrm{ref}}(s) \in SE(3)$, we also define the \textcolor{red}{strains}
\begin{equation*}
(\Lambda, \Gamma) = (K - K_{\mathrm{ref}}, W - W_{\mathrm{ref}}).
\end{equation*}
The simple case of a straight initial configuration along the $e_1$ axis, $g_{\mathrm{ref}}(s) = (I, s e_1)$, where $I$ is the identity matrix, leads to $\Lambda = K$ and $\Gamma = W - e_1$. One can see that $\Lambda$ measures the curvature (bending and torsion) and $\Gamma$ measures the difference between $d_1$ and $x^{\prime}$ (elongation and shear).}\\

\textcolor{blue}{Considering a hyperelastic material model with moderate strains, the Lagrangian density of the system can be written as
\begin{equation*}
\mathcal{L}(g,\dot{g},g^{\prime}) = \frac{1}{2}\left( \Omega^T \mathbb{J} \Omega + \rho A V^T V - \Lambda^T \mathbb{C}_1 \Lambda - \Gamma^T \mathbb{C}_2 \Gamma \right) - U_{\mathrm{ext}}(R, x)
\end{equation*}
where $\rho > 0$ is the linear density of the beam, $U_{\mathrm{ext}}: SE(3) \to \mathbb{R}$ is an external potential function and $\mathbb{J} = \rho \,\text{diag}([J_1, J_2, J_3])$ is the matrix of moments of inertia of the sections in the body frame. Assuming uniform cross sections and directors $d_2$ and $d_3$ coincident with the principal moments of area $I_2$ and $I_3$, one gets that $J_1 = \rho\,(I_2 + I_3)$, $J_2 = \rho I_2$ and $J_3 = \rho I_3$, and $\mathbb{C}_1 = \text{diag}([G (I_2 + I_3), E I_2, E I_3])$, $\mathbb{C}_2 = \text{diag}([E A, \kappa_2  G A, \kappa_3 G A])$, which are the matrices representing the corresponding stiffness parameters of the sections. $\textcolor{red}{\kappa_2}$ and $\textcolor{red}{\kappa_3}$ are possible shear correction factors.}\\

\subsection{\textcolor{blue}{Unit dual quaternion formulation}}
\textcolor{blue}{Working on $SE(3)$ is difficult. In \cite{Leitz2021} the authors propose the use of a constrained approach where the space of dual quaternions, $\widetilde{\mathbb{H}}$, which is a vector space, is considered as ambient manifold and the unit dual quaternions, $\widetilde{\mathbb{H}}_1$, as constraint submanifold since it is well-known that this latter space provides a double covering of $SE(3)$.}\\

\textcolor{blue}{The space of dual quaternions is defined by
\begin{equation*}
\widetilde{\mathbb{H}} := \left\lbrace
\tilde{q} = q_r + q_{d} \boldsymbol{\epsilon} \, \vert \, q_r, q_{d} \in \mathbb{H}, \boldsymbol{\epsilon}^2 = 0
\right\rbrace
\end{equation*}
where $\mathbf{\epsilon}$ is the so-called dual unit and
\begin{equation*}
\mathbb{H} := \left\lbrace
q = q_0 + q_1 \boldsymbol{i} + q_2 \boldsymbol{j} + q_3 \boldsymbol{k} \, \vert \, q_i \in \mathbb{R}, \boldsymbol{i}^2 = \boldsymbol{j}^2 = \boldsymbol{k}^2 = \boldsymbol{i j k} = -1
\right\rbrace
\end{equation*}
is the space of quaternions. Both of these are vector spaces, so working with them is quite simple.}\\

\textcolor{blue}{Similar to complex numbers, a conjugation operation is defined on the space of quaternions, namely, if $p = p_0 + p_1 \boldsymbol{i} + p_2 \boldsymbol{j} + p_3 \boldsymbol{k}$, then $\bar{p} = p_0 - q_1 \boldsymbol{i} - q_2 \boldsymbol{j} - q_3 \boldsymbol{k}$, and this operation is inherited by dual quaternions. This defines a norm on $\mathbb{H}$, $\Vert p \Vert = \sqrt{\bar{p} p}$ and lets us write the inverse of $p$ as $p^{-1} = \bar{p}/\Vert p \Vert^2$. This also defines a seminorm on $\widetilde{\mathbb{H}}$ by $\Vert \tilde{p} \Vert = \sqrt{\bar{\tilde{p}} \tilde{p}} = \Vert p_r \Vert + \frac{p_r^T p_\epsilon}{\Vert p_r \Vert} = \sqrt{p_r^T p_r} + \frac{p_r^T p_\epsilon}{\sqrt{p_r^T p_r}}$, where in the last equality we consider the quaternions $q_r, q_{\epsilon}$ as vectors in $\mathbb{R}^4$. The set of unit quaternions and unit dual quaternions are thus, $\mathbb{H}_1 := \left\lbrace
q \in \mathbb{H} \, \vert \, \Vert q \Vert = 1
\right\rbrace$ and $\widetilde{\mathbb{H}}_1 := \left\lbrace
\tilde{q} \in \widetilde{\mathbb{H}} \, \vert \, \Vert \tilde{q} \Vert = 1
\right\rbrace$ respectively. More explicitly, the latter implies
\begin{subequations}
\label{eq:unity_constraints}
\begin{align}
q_{0,r}^2 + q_{1,r}^2 + q_{2,r}^2 + q_{3,r}^2 &= 1\\
q_{0,r} q_{0,\epsilon} + q_{1,r} q_{1,\epsilon} + q_{2,r} q_{2,\epsilon} + q_{3,r} q_{3,\epsilon} &= 0
\end{align}
\end{subequations}}\\

\textcolor{blue}{As stated before, an element $\tilde{q} \in \widetilde{\mathbb{H}}_1$ can be put into correspondence with an element of $SE(3)$. In particular, we can parametrize $\tilde{q}$ by a rotation angle $\theta$ and two purely imaginary quaternions $n, x$, i.e. $n_0 = x_0 = 0$, with $\Vert n \Vert = 1$, representing a rotation axis and a three dimensional translation respectively. This way $q = \cos(\theta/2) + n \sin(\theta/2)$ and $\tilde{q} = q + \frac{1}{2} x q \boldsymbol{\epsilon}$. If $\tilde{q}(t,s) \in \widetilde{\mathbb{H}}_1$, then
\begin{equation*}
\widetilde{\Omega} := 2 \bar{\tilde{q}}\dot{\tilde{q}} = \Omega + V \boldsymbol{\epsilon}\,, \qquad \widetilde{K} := 2 \bar{\tilde{q}}\tilde{q}^{\prime} = K + W \boldsymbol{\epsilon}
\end{equation*}}

\textcolor{blue}{One can thus define an ambient Lagrangian in the dual quaternions,
\begin{equation}
	\label{eq:ambient_lagrangian_beam}
  \widetilde{\mathcal{L}}\left(\tilde{q},\dot{\tilde{q}},\tilde{q}^{\prime}\right) = 2 \widetilde{M}\left(\bar{\tilde{q}}\dot{\tilde{q}}, \bar{\tilde{q}}\dot{\tilde{q}}\right) - 2 \widetilde{C}\left( \bar{\tilde{q}}\tilde{q}^{\prime} - \bar{\tilde{q}}_{\mathrm{ref}}\tilde{q}_{\mathrm{ref}}^{\prime}, \bar{\tilde{q}}\tilde{q}^{\prime} - \bar{\tilde{q}}_{\mathrm{ref}}\tilde{q}_{\mathrm{ref}}^{\prime}\right) - \widetilde{U}(\tilde{q})
\end{equation}
where $\widetilde{M}(\tilde{q},\tilde{p}) = q_{r}^T \widetilde{\mathbb{J}} p_{r} + q_{\epsilon}^T \tilde{\rho} p_{\epsilon}$, $\widetilde{C}(\tilde{q},\tilde{p}) = q_{r}^T \widetilde{\mathbb{C}}_1 p_{r} + q_{\epsilon}^T \widetilde{\mathbb{C}}_2 p_{\epsilon}$, with $\widetilde{\mathbb{J}} = \text{diag}([\alpha_1,\mathbb{J}])$, $\tilde{\rho} = \text{diag}([\alpha_2,\rho A I])$, $\widetilde{\mathbb{C}}_1 = \text{diag}([\alpha_3,\mathbb{C}_1])$ and $\widetilde{\mathbb{C}}_2 = \text{diag}([\alpha_4,\mathbb{C}_2])$, and $\alpha_i \in \mathbb{R}$. These $\alpha$ can be chosen arbitrarily as they play no role in the dynamics once the unity constraints \eqref{eq:unity_constraints} are enforced.}

\subsection{\textcolor{blue}{Discrete Lagrangian}}
\begin{figure}%
    \centering%
    \includegraphics[scale=1.2]{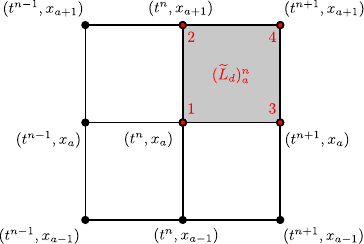}%
    \caption{Spacetime grid for the multisymplectic variational integrator \textcolor{blue}{with relative indices marked in red}.}
    \label{fig:spacetime}
\end{figure}

In order to discretize the beam, the spacetime $[0,T] \times [0, \ell]$ is discretized into a regular grid (see Figure \ref{fig:spacetime}) with constant space and time steps, $\Delta s$ and $\Delta t$ respectively.\\

\textcolor{blue}{We discretize the ambient Lagrangian density \eqref{eq:ambient_lagrangian_beam} applying the trapezoidal rule in both space and time
\begin{align*}
&\widetilde{L}_d(\tilde{q}_a^n, \tilde{q}_{a+1}^n, \tilde{q}_a^{n+1}, \tilde{q}_{a+1}^{n+1}) = \frac{1}{4 \Delta t \Delta s} \\
&\quad \times \left[ \widetilde{\mathcal{L}}\left( \tilde{q}_a^n, \frac{\tilde{q}_a^{n+1} - \tilde{q}_a^n}{\Delta t}, \frac{\tilde{q}_{a+1}^{n} - \tilde{q}_a^n}{\Delta s}\right) + \widetilde{\mathcal{L}}\left( \tilde{q}_{a+1}^n, \frac{\tilde{q}_{a+1}^{n+1} - \tilde{q}_{a+1}^n}{\Delta t}, \frac{\tilde{q}_{a+1}^{n} - \tilde{q}_a^n}{\Delta s}\right)\right.\\
&\quad\left. + \,\widetilde{\mathcal{L}}\left( \tilde{q}_a^{n+1}, \frac{\tilde{q}_a^{n+1} - \tilde{q}_a^n}{\Delta t}, \frac{\tilde{q}_{a+1}^{n+1} - \tilde{q}_a^{n+1}}{\Delta s}\right) + \widetilde{\mathcal{L}}\left( \tilde{q}_{a+1}^{n+1}, \frac{\tilde{q}_{a+1}^{n+1} - \tilde{q}_{a+1}^n}{\Delta t}, \frac{\tilde{q}_{a+1}^{n+1} - \tilde{q}_a^{n+1}}{\Delta s}\right)\right]
\end{align*}
and introduce the notation $(\widetilde{L}_d)_a^n := \widetilde{L}_d(\tilde{q}_a^n, \tilde{q}_{a+1}^n, \tilde{q}_a^{n+1}, \tilde{q}_{a+1}^{n+1})$ to simplify the formulas.}

\textcolor{blue}{As derived in \cite{Leitz2021}, the discrete constrained \textrm{\textsc{Euler}}-\textrm{\textsc{Lagrange}} field equations are derived via a discrete variational principle in space and time and subsequent rearrangement of terms in space index $a$ and time index $n$. As shown there, a natural choice of null space matrix is
\begin{equation*}
\widetilde{P}(\tilde{q}) = \left[
\begin{array}{cc}
P(q_r) & 0\\
P(q_{\epsilon}) & P(q_r)\\
\end{array}
\right]\,, \qquad
P(q) = \frac{1}{2} \left[
\begin{array}{rrr}
-q_1 & -q_2 & -q_3\\
 q_0 & -q_3 &  q_2\\
 q_3 &  q_0 & -q_1\\
-q_2 &  q_1 &  q_0
\end{array}
\right]\,.
\end{equation*}}

\begin{figure}
    \centering
    \includegraphics[scale=0.4]{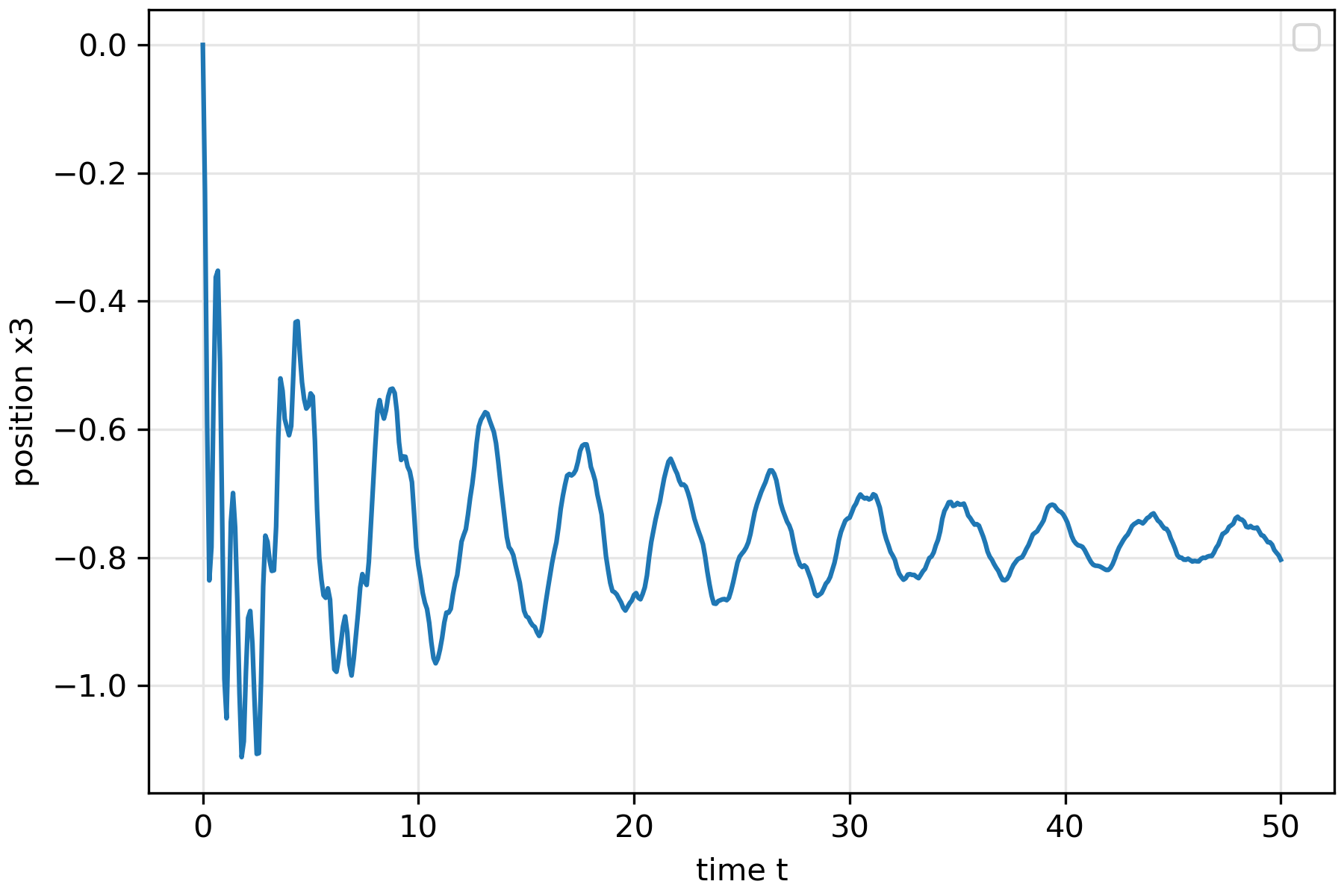}
    \caption{Viscous damping of a cantilever beam.}
    \label{fig:diss}
\end{figure}

\textcolor{blue}{The forced version of these equations results from the application of the discrete  \textrm{\textsc{Lagrange}}-\textrm{\textsc{d'Alembert}} principle, similar to \eqref{eqn:LAprinc},
\begin{equation*}
\sum_{a} \sum_{n} \left[ \delta (\widetilde{L}_d)_{a}^n + (f_d^1)_{a}^{n} \delta \tilde{q}_{a}^n  + (f_d^2)_{a}^{n} \delta \tilde{q}_{a+1}^n + (f_d^3)_{a}^{n} \delta \tilde{q}_{a}^{n+1} + (f_d^4)_{a}^{n} \delta \tilde{q}_{a+1}^{n+1} \right] = 0\,.
\end{equation*}
with $(f_d^i)_{a}^n := f_d^i(\tilde{q}_a^n, \tilde{q}_{a+1}^n, \tilde{q}_a^{n+1}, \tilde{q}_{a+1}^{n+1},u_a^n)$ denoting all external and control forces, and $i = 1,...,4$ coinciding with the corresponding relative node on which they are applied, as in Figure~\ref{fig:spacetime}. This leads to DEL with a force contribution from each adjacent spacetime rectangle sharing the node under consideration:
\begin{equation}
    \begin{aligned}
    \widetilde{P}(\tilde{q}_a^n)^T \,\Big[ D_1 (\widetilde{L}_d)_a^n + D_2 (\widetilde{L}_d)_{a-1}^n + D_3 (\widetilde{L}_d)_a^{n-1} + D_4 (\widetilde{L}_d)_{a-1}^{n-1} &\\
    + (f_d^1)_{a}^{n} + (f_d^2)_{a-1}^{n} + (f_d^3)_{a}^{n-1} + (f_d^4)_{a-1}^{n-1} \Big] &= 0.
    \end{aligned}
\end{equation}}
Suitable boundary conditions in space and time as well at the spacetime corners are directly derived via the discrete variational principle.\\
\textrm{\textsc{Kelvin}}-\textrm{\textsc{Voigt}} type viscous damping is included as external forces that are proportional to the discrete approximation of the strain rate \cite{Linn2013} with bulk viscosity \textcolor{blue}{$\zeta$ and shear viscosity $\eta$. In the moderate strain regime these result in a damping matrix $\widetilde{\mathbb{D}} = \mathrm{diag}([0,\eta (I_2 + I_3), \chi I_2, \chi I_3, 0, \chi A, \eta A, \eta A])$, where $\chi = \zeta (3 - E/G)^2 + \eta (E/G)^2/3$ is the extensional viscosity. The corresponding discrete force is
\begin{align*}
\mathrm{L}_{(\tilde{q}^{n}_a)}^T {(f^{\mathrm{KV, 1}}_d)^n_a} = \mathrm{L}_{(\tilde{q}^{n+1}_a)}^T {(f^{\mathrm{KV, 3}}_d)^n_a} &= \frac{\Delta t \Delta s}{4}\,
\widetilde{\mathbb{D}}\left(\frac{\widetilde{K}^{n+1}_a - \widetilde{K}^n_a}{\Delta t}\right),\\
\mathrm{L}_{(\tilde{q}^{n}_{a+1})}^T {(f^{\mathrm{KV, 2}}_d)^n_a} = \mathrm{L}_{(\tilde{q}^{n+1}_{a+1})}^T {(f^{\mathrm{KV, 4}}_d)^n_{a+1}} &= \frac{\Delta t \Delta s\,}{4}
\widetilde{\mathbb{D}}\left(\frac{\widetilde{K}^{n+1}_{a+1} - \widetilde{K}^n_{a+1}}{\Delta t}\right),
\end{align*}
where by $\mathrm{L}_{\tilde{q}}^T$, we denote the transposed of the dual quaternion left multiplication operation by $\tilde{q}$, and $\widetilde{K}^{n}_a = \widetilde{K}(\tilde{q}^{n}_a,(\tilde{q}')^{n}_a)$, with $(\tilde{q}')^n_a = (\tilde{q}')^n_{a+1} = (\tilde{q}^n_{a+1} - \tilde{q}^n_{a})/\Delta s$ and $(\tilde{q}')^{n+1}_a = (\tilde{q}')^{n+1}_{a+1} = (\tilde{q}^{n+1}_{a+1} - \tilde{q}^{n+1}_{a})/\Delta s$.}
Figure \ref{fig:diss} shows the position of the tip of a cantilever beam with fixed-free boundary conditions that is initially straight under gravity. The strain-rate proportional damping leads to reduced high frequency oscillations. 
\\
\subsection{The Discrete Adjoint Method in Spacetime}
The discrete adjoint method for the geometrically exact beam considers the configuration variables as well as the adjoint variables in space and time to derive the discrete adjoint equations in space and time. Single shooting in time while simultaneously solving the equations in space is used for the solution of the optimal control problem. The \textrm{\textsc{Barzilei}}-\textrm{\textsc{Borwein}} gradient method \cite{Barzilai1988,Fletcher2001} is used for the update. Here, a pendulum-like beam \textcolor{blue}{subject to} gravity and fixed-free translation and free-free rotation boundary conditions is considered with a torque \textcolor{blue}{$u$} applied at the fixed end \textcolor{blue}{as discrete redundant control forces
\begin{equation*}
\mathrm{L}_{(\tilde{q}^{n}_0)}^T (f^{\mathrm{C, 1}}_d)_0^n = \mathrm{L}_{(\tilde{q}^{n+1}_0)}^T (f^{\mathrm{C, 3}}_d)_0^n = 2 \Delta t \Delta s \, u_0^n \, \boldsymbol{k}\,.
\end{equation*}
Since our control is only applied at the boundary, this is a boundary control problem for a PDE.}

\textcolor{blue}{The desired configuration is the upright rotated position of the beam, specified for each node in space. The final position considered is undeformed. The desired maneuver is from the lower position to the upright position in such a way that the inertial terms cancel the strains in the end configuration. As the system is heavily underactuated, the chosen input does not allow us to control the motion in the axial direction and does not lead to a stationary upright position. Hence, no end momentum is imposed. Nonetheless, the control task should demonstrate the presented method in an academic example that resembles the previous pendulum examples sufficiently.}

\textcolor{blue}{
Our optimal control problem is of the form
\begin{subequations}\label{eqn:field_OCP}
    \begin{equation}
        \underset{\tilde{q}_d, u_d}{\min}\, J_d(\tilde{q}_d, u_d) = \sum_{a = 0}^A [\tilde{q}_a^N - (\tilde{q}_a^N)_{*}]^T S_q [\tilde{q}_a^N - (\tilde{q}_a^N)_{*}] + \sum_{n=0}^{N-1} \textcolor{blue}{\frac{1}{2}} (u_0^n)^T R (u_0^n)
    \end{equation}
subject to:
    \begin{align}
        \tilde{q}_a^0 &= (\tilde{q}_a^0)_*, \hspace{6.4cm}\text{for}~a=0,~...,~N \\
        u_0^0 &= (u_0^0)_*,\\
        0 &= \textcolor{red}{\Delta s (p_a^0)_* +} \widetilde{P}^T(\tilde{q}_a^0) \left[ D_1 (\widetilde{L}_d)_a^0 + D_2 (\widetilde{L}_d)_{a-1}^0 + (f_{d}^{1})_a^0 + (f_d^2)_{a-1}^{0}\right],\\
        &\hspace{8.1cm}\text{for}~a=1,~...,~A-1 \nonumber\\
        0 &= \widetilde{P}^T(\tilde{q}_0^n) \left[ D_1 (\widetilde{L}_d)_0^n + D_3 (\widetilde{L}_d)_0^{n-1} + (f_{d}^{1})_0^n + (f_d^3)_{0}^{n-1}\right],\\
        &\hspace{8.1cm}\text{for}~n=1,~...,~N-1 \nonumber\\
        0 &= \widetilde{P}^T(\tilde{q}_a^n) \left[ D_1 (\widetilde{L}_d)_a^n + D_2 (\widetilde{L}_d)_{a-1}^n + D_3 (\widetilde{L}_d)_a^{n-1} + D_4 (\widetilde{L}_d)_{a-1}^{n-1}\right.\\
          &+\left. (f_d^1)_{a}^{n} + (f_d^2)_{a-1}^{n} + (f_d^3)_{a}^{n-1} + (f_d^4)_{a-1}^{n-1}\right],\nonumber\\
          &\hspace{4.8cm}\text{for}~a=1,~...,~A-1,\,\text{for}~n=1,~...,~N-1 \nonumber
    \end{align}
\end{subequations}
where $(\tilde{q}_a^0)_{*}$, $(u_0^0)_{*}$, \textcolor{red}{$(p_a^0)_*$}, are given initial discrete values and $(\tilde{q}_a^N)_{*}$ denotes the discretised desired end configuration. \textcolor{red}{$(p_a^0)_*$ are discrete initial temporal momenta, canonically associated to discretised linear and angular velocities, which in our example are set to $0$.}}

The adjoint equations are obtained similar to the constrained temporal case by applying discrete variational calculus and nodal reparametrization, but now in space and time. \textcolor{blue}{However, the resulting equations are quite long, and so, will not be reproduced here in their entirety. For instance, the equations obtained by taking variations of the inputs at the fixed boundary $a=0$, are
\begin{subequations}\label{eqn:BeamInp}
    \begin{align}
        &0= (\lambda_0^{n})^{T}  \widetilde{P}^T(\tilde{q}_0^n) (D_5 f_d^1)_0^n  +  (\lambda_0^{n+1})^{T} \widetilde{P}^T(\tilde{q}_0^{n+1}) (D_5 f_d^2)_{0}^{n}\,,\\
        &\hspace{8.65cm} \text{for}~n=1,~...,~N-2 \nonumber\\
        &0= (\lambda_0^{N-1})^{T}  \widetilde{P}^T(\tilde{q}_0^{N-1}) (D_5 f_d^1)_{0}^{N-1}.
    \end{align}
\end{subequations}
These are used to update the torque. If instead of boundary control we had controls over the bulk, these equations would generalise to all nodes as:
\begin{equation*}
    \begin{aligned}
    0&=(\lambda_{a}^{n})^{T} \widetilde{P}^T(\tilde{q}_a^{n}) (D_5 f_d^1)_{a}^{n} + (\lambda_{a+1}^{n})^{T} \widetilde{P}^T(\tilde{q}_{a+1}^{n}) (D_5 f_d^2)_{a}^{n}\\
    &+ (\lambda_{a}^{n+1})^{T} \widetilde{P}^T(\tilde{q}_a^{n+1}) (D_5 f_d^3)_{a}^{n} + (\lambda_{a+1}^{n+1})^{T} \widetilde{P}^T(\tilde{q}_{a+1}^{n+1}) (D_5 f_d^4)_{a}^{n},\\
    &\hspace{5.5cm}  \mathrm{for}~a=1,~...,~A-2,~n=1,~...,~N-2\\
    0&=(\lambda_{A-1}^{n})^{T} \widetilde{P}^T(\tilde{q}_{A-1}^{n}) (D_5 f_d^1)_{A-1}^{n} + (\lambda_{A-1}^{n+1})^{T} \widetilde{P}^T(\tilde{q}_{A-1}^{n+1}) (D_5 f_d^3)_{A-1}^{n},\\
    &\hspace{8.35cm} \text{for}~n=1,~...,~N-2\\
    0&=(\lambda_{a}^{N-1})^{T} \widetilde{P}^T(\tilde{q}_a^{N-1}) (D_5 f_d^1)_{a}^{N-1} + (\lambda_{a+1}^{N-1})^{T} \widetilde{P}^T(\tilde{q}_{a+1}^{N-1}) (D_5 f_d^2)_{a}^{N-1},\\
    &\hspace{8.45cm} \text{for}~a=1,~...,~A-2\\
    0&=(\lambda_{A-1}^{N-1})^{T} \widetilde{P}^T(\tilde{q}_{A-1}^{N-1}) (D_5 f_d^1)_{A-1}^{N-1}.
    \end{aligned}
\end{equation*}
}

\subsection{\textcolor{blue}{Fairly rigid} Beam}
\begin{figure}
\centering
	\begin{subfigure}{0.45\textwidth}
	  \centering
	\centering\resizebox{0.9\textwidth}{!}{\includegraphics{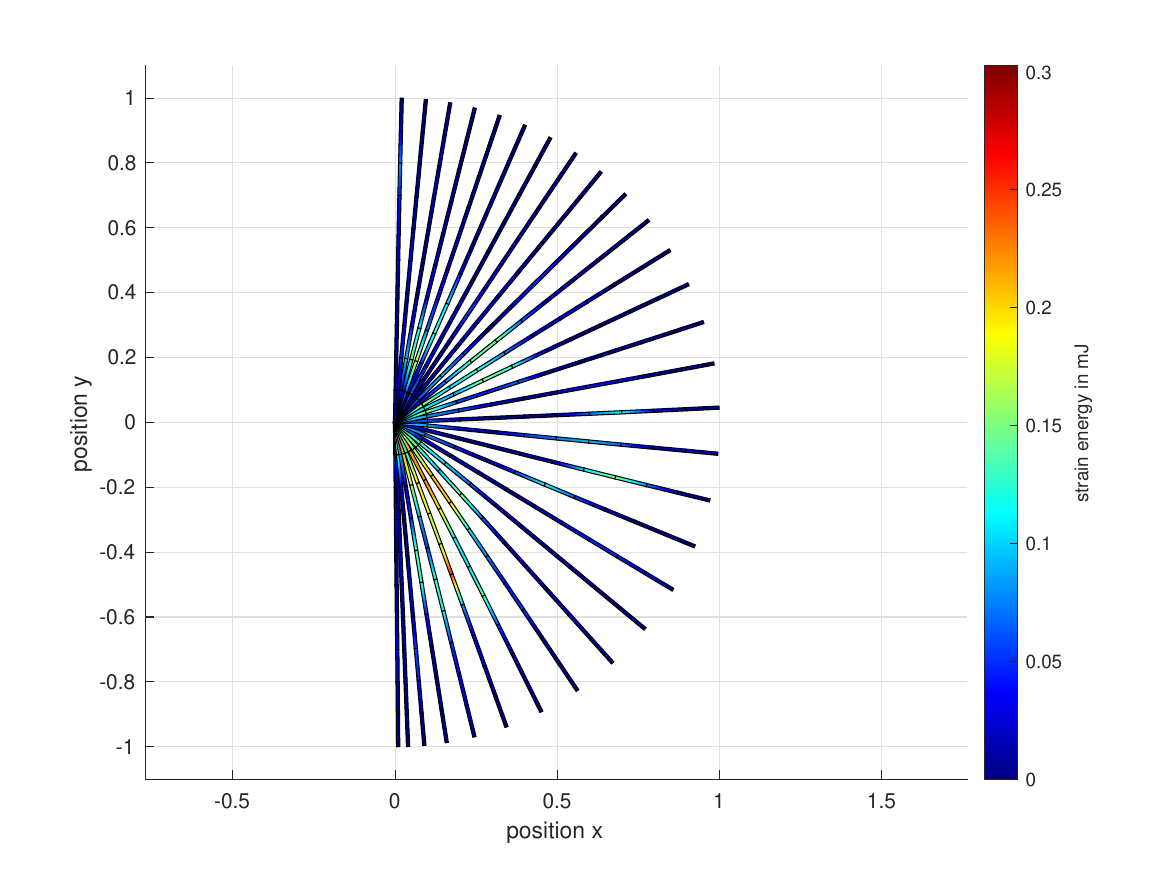}}
\caption{Deformation of a \textcolor{blue}{fairly rigid} beam during the \textcolor{blue}{upward motion}.}\label{subfig:steel_motion}
	\end{subfigure}%
 \quad
	\begin{subfigure}{0.45\textwidth}
	  \centering
	\centering\resizebox{0.9\textwidth}{!}{\input{steel/energy.tex}}
\caption{Energy of the beam for the optimized \textcolor{blue}{motion} over time $t$.}\label{subfig:steel_energy}
	\end{subfigure}
	\caption{\textcolor{blue}{Motion} of a \textcolor{blue}{fairly rigid} beam and its energy.}
\label{fig:steelBeam}
\end{figure}

The \textcolor{blue}{fairly rigid} beam demonstrates the sequential optimization of the beam dynamics with objective minimization of the control effort. The simulation of the beam dynamics uses $A=10$ nodes in space and $N=3000$ nodes in time. The beam has a length of $L=1$. The simulation duration is $T=1$. The resulting time step is $h = \frac{1}{3000}$ in time and the step size in space is $\Delta s = \frac{1}{10}$. A constant initial guess of $u^0=1500$ is used. The beam has a square cross-section of $\textcolor{blue}{A_{\mathrm{cross}}}=0.01$ with a side length of $\textcolor{blue}{l_s}=0.1$. The \textcolor{blue}{chosen} weighting for the end term is $S_q=10^8$, and $R = 10^{-2}$ for the input\footnote{\textcolor{blue}{These values have been chosen to provide similar magnitudes to the different terms of the discrete objective. Notice that $S_q$ affects only a single time step, multiplying terms with values around $\pi$ and $0$. $R$ appears in a sum containing the $3000$ time steps with input values between $2500$ and $0$.}}. The material of the beam is \textcolor{blue}{fairly rigid} with a \textrm{\textsc{Young's}} modulus of \textcolor{blue}{$E=210000$} and a \textrm{\textsc{Poisson}} ratio of $\nu=0.3$. The mass density is $\rho = 7.85$. The beam is damped with $\eta = 1\cdot10^{-1}$ and $\zeta=1\cdot10^{-2}$.\\
Figure \ref{subfig:steel_motion} shows snapshot of the \textcolor{blue}{motion} of the beam. Figure \ref{subfig:steel_energy} shows the total energy \textcolor{blue}{$H$} as well as all its contributions over time. The deformation energy is the difference between the \textcolor{blue}{total} potential energy \textcolor{blue}{of the system} $\textcolor{blue}{U}$ and the gravitational potential energy $\textcolor{blue}{U}_{grav}$. \textcolor{blue}{The main contribution to the kinetic energy $T$ is due to translation. At the end of the simulation, the kinetic energy reduces due to the input weight.} The optimized input is depicted in Figure \ref{subfig:steel_inp}, it decreases to zero at the end of the simulation time. The optimized quantities, the distance of the beam to the desired end configuration as well as the control effort are depicted in Figures \ref{subfig:steel_dist} and \ref{subfig:steel_ceff}, respectively. \textcolor{blue}{The gradient depicted in Figure \ref{subfig:steel_grad} shows heavy oscillations. }
\begin{figure}
\centering
	\begin{subfigure}{0.45\textwidth}
	  \centering
	\centering\resizebox{0.9\textwidth}{!}{\input{steel/input.tex}}
	\caption{Input $u_n$ over time $t$.}\label{subfig:steel_inp}
	\end{subfigure}%
 \quad
	\begin{subfigure}{0.45\textwidth}
	  \centering
	\centering\resizebox{0.9\textwidth}{!}{
%
%
\definecolor{mycolor1}{rgb}{0.00000,0.44700,0.74100}%
\begin{tikzpicture}

\begin{axis}[%
xmin=0,
xmax=67,
xlabel style={font=\color{white!15!black}},
xlabel={iterations},
ymode=log,
ymin=0.01,
ymax=1.94652417295488,
yminorticks=true,
ylabel style={font=\color{white!15!black}},
ylabel={distance $\sum_{a=0}^{A}\left\Vert q_a^J-q_a^{J~\mathrm{pre}}\right\Vert$},
axis background/.style={fill=white},
xmajorgrids,
ymajorgrids,
yminorgrids
]
\addplot [color=mycolor1, forget plot]
  table[row sep=crcr]{%
1	1.94652417295488\\
2	0.47194864902709\\
3	0.418243415005795\\
4	0.0908473793535848\\
5	0.0545472156460274\\
6	0.0575335731204717\\
7	0.0573608931292679\\
8	0.0249708239574426\\
9	0.617142779768869\\
10	0.0673177445769132\\
11	0.0496426646234501\\
12	0.0518474185390906\\
13	0.0517269466513231\\
14	0.0471844804132456\\
15	0.281623986696335\\
16	0.0922449176611592\\
17	0.0464192016728935\\
18	0.0486120681674409\\
19	0.04855989489653\\
20	0.0419885854641033\\
21	0.351765277824303\\
22	0.0732916707064649\\
23	0.0444259674276859\\
24	0.046134690632927\\
25	0.0460984967572148\\
26	0.0439309416577116\\
27	0.214203966872941\\
28	0.0845774811804788\\
29	0.0429418461579549\\
30	0.0442812285352115\\
31	0.0442610152098727\\
32	0.0430368801042738\\
33	0.164223265074223\\
34	0.0824663118134535\\
35	0.0418933395966926\\
36	0.042813101546689\\
37	0.0428018472269505\\
38	0.0422590306962623\\
39	0.106181773321817\\
40	0.0817288940992125\\
41	0.0409417001063698\\
42	0.0414783282996558\\
43	0.0414725886696553\\
44	0.0411755486337122\\
45	0.0698551091326523\\
46	0.0888970953545654\\
47	0.038765980987256\\
48	0.0393435023574716\\
49	0.0393361371492363\\
50	0.0391118043498372\\
51	0.0448908720688985\\
52	0.0595805694581904\\
53	0.0378269473518799\\
54	0.0385072322248824\\
55	0.0384977760909051\\
56	0.0383116780304707\\
57	0.040543698596184\\
58	0.0478823100814626\\
59	0.037280389140301\\
60	0.0381082947138115\\
61	0.0380907961359919\\
62	0.0379754332555552\\
63	0.0395376717555894\\
64	0.0455748923274178\\
65	0.0368825064456238\\
66	0.03783863978922\\
67	0.0378140022591684\\
};
\end{axis}

\end{tikzpicture}%
}
	\caption{Distance to end configuration over optimization iterations.}\label{subfig:steel_dist}
	\end{subfigure}%
    \bigbreak
	\begin{subfigure}{0.45\textwidth}
	  \centering
	\centering\resizebox{0.9\textwidth}{!}{
%
%
\definecolor{mycolor1}{rgb}{0.00000,0.44700,0.74100}%
\begin{tikzpicture}

\begin{axis}[%
xmin=0,
xmax=67,
xlabel style={font=\color{white!15!black}},
xlabel={iterations},
ymin=82000,
ymax=100000,
ylabel style={font=\color{white!15!black}},
ylabel={control effort $\left\Vert u_n \right\Vert$},
axis background/.style={fill=white},
axis x line*=bottom,
axis y line*=left,
xmajorgrids,
ymajorgrids
]
\addplot [color=mycolor1, forget plot]
  table[row sep=crcr]{%
1	82144.6894205584\\
2	99832.4111894806\\
3	94298.1576216969\\
4	96327.0462979216\\
5	96483.7940866504\\
6	96401.0761822805\\
7	96351.1652554341\\
8	92449.4241941015\\
9	88260.1193990647\\
10	91658.8487646145\\
11	91740.9176384077\\
12	91706.7686070311\\
13	91693.0214487511\\
14	91020.0806474341\\
15	88774.0537488046\\
16	89874.1720112965\\
17	90129.6345485152\\
18	90109.6887565905\\
19	90104.5295471412\\
20	89643.1910532888\\
21	87712.4227860623\\
22	89288.4296903058\\
23	89446.6252311073\\
24	89435.3478978672\\
25	89434.1957368996\\
26	89354.3058530591\\
27	88373.5059379032\\
28	89080.2346063308\\
29	89304.8567804217\\
30	89297.7577033969\\
31	89297.9518963133\\
32	89320.2749588469\\
33	88763.1973521985\\
34	89197.0612393003\\
35	89411.1801282769\\
36	89406.9090946315\\
37	89407.4091910795\\
38	89450.4866384735\\
39	89271.9787371309\\
40	89399.0225553062\\
41	89610.3393166848\\
42	89608.1104387903\\
43	89608.5446600911\\
44	89647.7699758446\\
45	89681.8644216082\\
46	89583.9543978415\\
47	89840.7651568082\\
48	89838.5637914228\\
49	89839.1079959783\\
50	89866.927867455\\
51	89915.7693390401\\
52	89839.2948872768\\
53	89950.5907939311\\
54	89947.9508549858\\
55	89948.5668496666\\
56	89969.9023458384\\
57	89996.5909342412\\
58	89956.962472248\\
59	90011.5400318164\\
60	90008.1934760463\\
61	90008.9009689656\\
62	90019.408289519\\
63	90041.6631401329\\
64	90008.2061081105\\
65	90053.0933751958\\
66	90049.1283392761\\
67	90049.9077125067\\
};
\end{axis}

\end{tikzpicture}%
}
	\caption{Control effort $\left\Vert u_n \right\Vert$ over optimization iterations.}\label{subfig:steel_ceff}
	\end{subfigure}
 \quad
     \begin{subfigure}{0.45\textwidth}
	  \centering
	\centering\resizebox{0.9\textwidth}{!}{\input{steel/grad.tex}}
	\caption{\textcolor{blue}{Gradient of the augmented objective with respect to the input $\frac{\partial \tilde{J}}{\partial u_n}$ over time $t$.}}\label{subfig:steel_grad}
	\end{subfigure}%
	\caption{Optimization results of beam \textcolor{blue}{dynamics} using single shooting for a \textcolor{blue}{fairly rigid} beam.}
	\label{fig:opt_steel}
\end{figure}

\subsection{\textcolor{blue}{Very flexible} Beam}
A \textcolor{blue}{very flexible} beam demonstrates the sequential optimization for more flexible beams that show larger deformations and are therefore harder to control. The simulation of the beam and adjoint dynamics uses $A=5$ nodes in space for a length of $L=1$. This results in a space step width of $\Delta s = \frac{1}{5}$. The simulation time is $T=0.5$ using $N=600$ node in time and a time step of $h = \frac{1}{1200}$. The initial guess for the input is $u^0 = 50$ for all time intervals. The beam has a square cross-section of $\textcolor{blue}{A_{\mathrm{cross}}}=0.0025$ with a side length of $l_s=0.05$. The \textrm{\textsc{Young's}} modulus is \textcolor{blue}{$E=50000$} and the mass density $\rho = 1000$. The \textrm{\textsc{Poisson}} ratio is $\nu=0.35$. \textrm{\textsc{Kelvin}}-\textrm{\textsc{Voigt}} type damping is applied with $\eta = 1\cdot10^{-1}$ and $\zeta=1\cdot10^{-2}$.\\
The weighting for the end configuration is $S_q=10^2$. \textcolor{blue}{For this numerical experiment, the input weight was set to $R=0$, since the chosen end configuration gets increasingly harder to reach for more flexible beams.}\\
\begin{figure}
\centering
	\begin{subfigure}{0.45\textwidth}
	  \centering
	\centering\resizebox{0.9\textwidth}{!}{\includegraphics{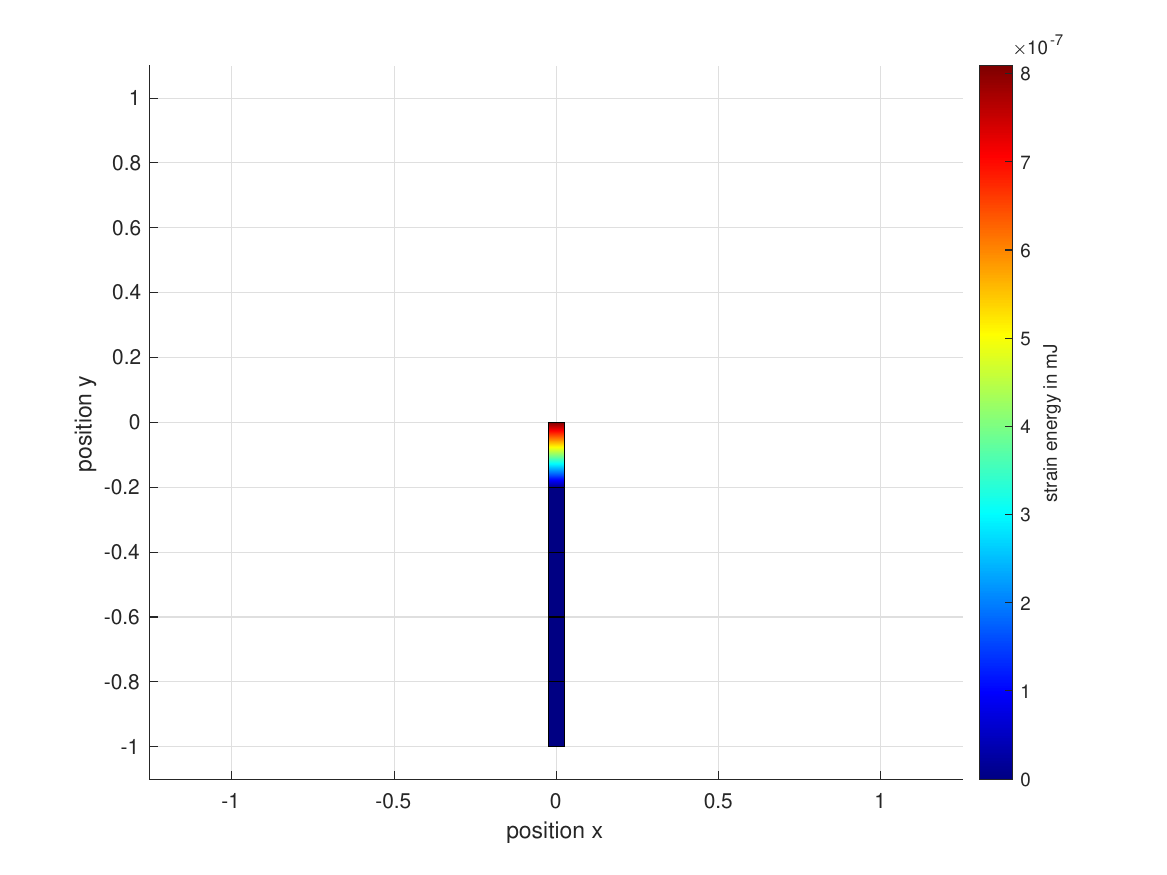}}
	\caption{Initial configuration of the beam by example of a \textcolor{blue}{very flexible} beam.}\label{subfig:iniConfig}
	\end{subfigure}
 \quad
	\begin{subfigure}{0.45\textwidth}
	  \centering
	\centering\resizebox{0.9\textwidth}{!}{\includegraphics{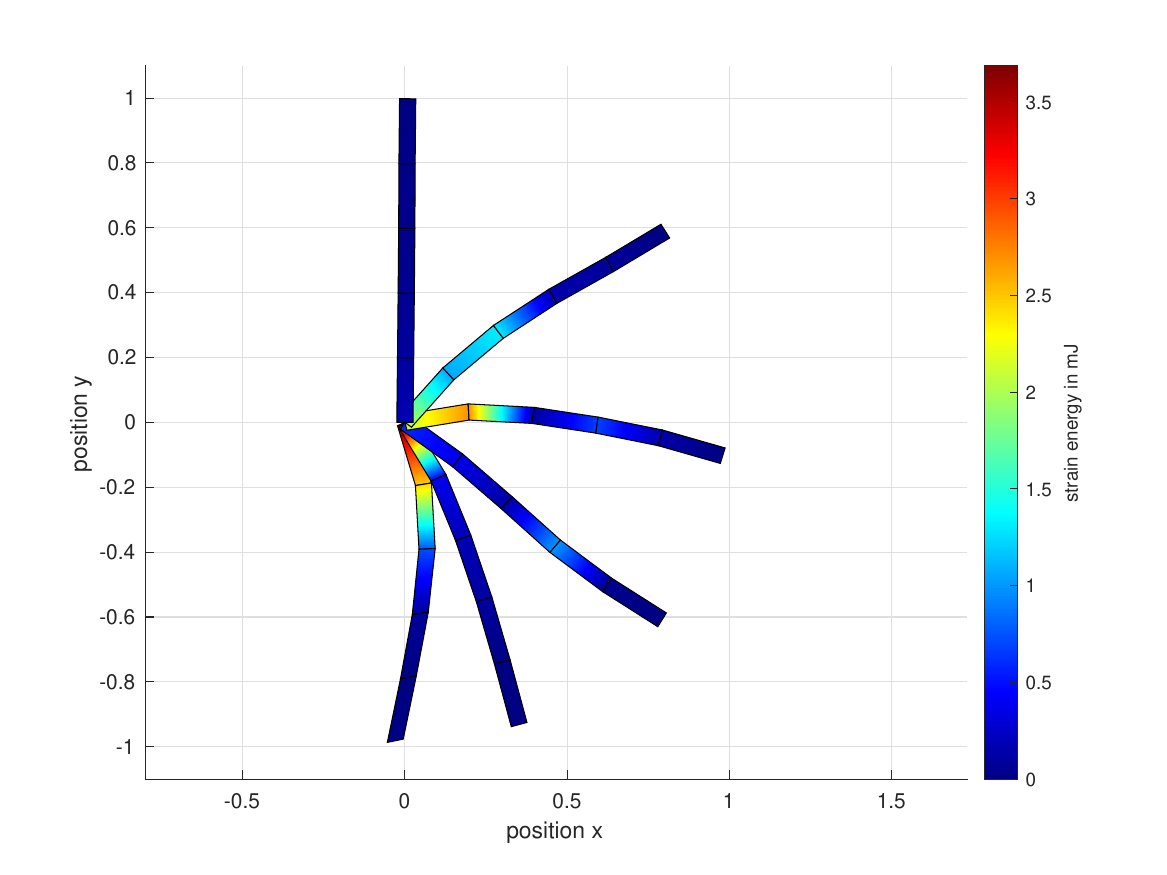}}
	\caption{Deformation of a \textcolor{blue}{very flexible} beam during \textcolor{blue}{its motion}.}\label{subfig:rubberBeam}
	\end{subfigure}%
\caption{Snapshot of a \textcolor{blue}{very flexible} beam at $n=0$ as well as during the \textcolor{blue}{ upward movement} of a \textcolor{blue}{very flexible} beam.}\label{fig:Beam}
\end{figure}
The optimization results are depicted in Figure \ref{fig:Beam}. The input in Figure \ref{subfig:rubber_inp} is increased compared to the initial guess. In addition, oscillations are present. 
The gradient depicted in Figure \ref{subfig:rubber_grad} shows oscillations with high frequency that are likely caused by the dynamics of the beam in normal direction as these deformations are of much higher frequency than bending deformations due to the difference in stiffness.
The objective is depicted in \ref{subfig:rubber_obj}. The largest decrease happens at the start of the optimization. Figure \ref{subfig:rubber_energy} depicts the total energy and its parts. During the optimization, mainly the translational part of the kinetic energy increases as well as the potential energy due to the gravitation.

\begin{figure}
\centering
	\begin{subfigure}{0.45\textwidth}
	  \centering
	\centering\resizebox{0.9\textwidth}{!}{\input{rubber/input.tex}}
	\caption{Input $u_n$ over time $t$.}\label{subfig:rubber_inp}
	\end{subfigure}
 \quad
	\begin{subfigure}{0.45\textwidth}
	  \centering
	\centering\resizebox{0.9\textwidth}{!}{\input{rubber/grad.tex}}
	\caption{Gradient of the augmented objective with respect to the input $\frac{\partial \tilde{J}}{\partial u_n}$ over time $t$.}\label{subfig:rubber_grad}
	\end{subfigure}
	\bigbreak
	\begin{subfigure}{0.45\textwidth}
	  \centering
	\centering\resizebox{0.9\textwidth}{!}{
%
%
\definecolor{mycolor1}{rgb}{0.00000,0.44700,0.74100}%
\begin{tikzpicture}

\begin{axis}[%
xmin=0,
xmax=113,
xlabel style={font=\color{white!15!black}},
xlabel={iterations},
ymode=log,
ymin=0.0001,
ymax=100,
yminorticks=true,
ylabel style={font=\color{white!15!black}},
ylabel={objective $J$},
axis background/.style={fill=white},
xmajorgrids,
ymajorgrids,
yminorgrids
]
\addplot [color=mycolor1, forget plot]
  table[row sep=crcr]{%
1	69.1512963220325\\
2	68.9642900252709\\
3	2.63253594188961\\
4	6.31633278005024\\
5	0.185198082363964\\
6	0.144269476088675\\
7	0.132993318674881\\
8	0.129125822245928\\
9	0.128365334177448\\
10	0.10347126908332\\
11	0.0643591984063604\\
12	0.0594019495487053\\
13	0.065638253768825\\
14	0.0740892968221978\\
15	0.0644511591982581\\
16	0.0575849482479016\\
17	0.055891911767375\\
18	0.0733898241436479\\
19	0.114184271061519\\
20	0.117216859326808\\
21	0.0982070880841693\\
22	0.0942355317995154\\
23	0.0914680206185536\\
24	0.0680801441982376\\
25	0.0459198167634155\\
26	0.0416336176490179\\
27	0.0457685301753765\\
28	0.0457684488957462\\
29	0.0414330324297787\\
30	0.0383148649846819\\
31	0.0557061540248943\\
32	0.0847158466184002\\
33	0.102951277440232\\
34	0.0843781962704103\\
35	0.0807881596510745\\
36	0.0785574076252015\\
37	0.0542471287402577\\
38	0.0364780639266718\\
39	0.0432053577109462\\
40	0.0254352849344982\\
41	0.0235967359725451\\
42	0.021267361374133\\
43	0.0122110353477927\\
44	0.0119221590559776\\
45	0.00883117600509565\\
46	0.00573068461568802\\
48	0.00553677940369453\\
49	0.00832446470665351\\
50	0.0150888239272746\\
51	0.0146981235931628\\
52	0.0143931496528157\\
53	0.0112679543389755\\
54	0.00190735632480398\\
55	0.00635922367667322\\
56	0.00130296206166396\\
57	0.00106797081302866\\
58	0.00105498369728482\\
59	0.000973799949529636\\
60	0.00264970177913185\\
61	0.000934127426032261\\
63	0.000935596719232526\\
64	0.000705905121381076\\
65	0.169362049459898\\
66	0.0116624189254031\\
67	0.00152054002504643\\
68	0.000912574130256316\\
70	0.000900640304402795\\
71	0.000864854758656848\\
72	0.000848061277969388\\
73	0.000971252428733854\\
74	0.000812278872192047\\
75	0.000816897954498224\\
76	0.000815959579420168\\
77	0.000821643449234542\\
78	0.000896250861287006\\
79	0.000861964376729332\\
80	0.000869343290237101\\
82	0.000862897671373142\\
83	0.000872382048557112\\
84	0.000889418873426539\\
85	0.000889391061748968\\
86	0.000895771680102669\\
88	0.000894636605082941\\
90	0.000911415383255216\\
93	0.000912761712240056\\
98	0.000904830667427784\\
100	0.000893939310862748\\
101	0.000899303318636238\\
102	0.000898331652349361\\
104	0.000884247258491037\\
105	0.000889449438677745\\
106	0.000879635612115864\\
108	0.000879124547241383\\
110	0.000870377612389859\\
111	0.000876575092639858\\
112	0.000866081884125423\\
113	0.000868438263190813\\
};
\end{axis}

\end{tikzpicture}%
}
	\caption{Objective $J$ over optimization iterations.}\label{subfig:rubber_obj}
	\end{subfigure}
 \quad
	\begin{subfigure}{0.45\textwidth}
	  \centering
	\centering\resizebox{0.9\textwidth}{!}{\input{rubber/energy.tex}}
	\caption{Energy of the beam for the optimized \textcolor{blue}{motion over time }$t$.}\label{subfig:rubber_energy}
	\end{subfigure}%
	\caption{Optimization results of beam \textcolor{blue}{dynamics of a very flexible beam}.}
	\label{fig:opt_rubber}
\end{figure}
\section{Summary}%
The discrete adjoint method for variational integration of (constrained) ODEs is derived, and its convergence properties are demonstrated with the help of numerical examples. Quadratic convergence results of the configuration variables as well as for the adjoint variables based on simulations of a mathematical pendulum are observed. The discrete adjoint method is also applied to the multisymplectic Galerkin Lie group integrator for geometrically exact beam dynamics, in particular to the optimal control of the \textcolor{blue}{upward motion} of a pendulum-like beam. The discrete adjoint method directly derives fitting equations at the boundary based on the discretization chosen for the variational integrator. The discrete adjoint method for constrained systems with null space projection and nodal reparametrization also directly results in the null space projection of the discrete adjoint equations. The properties of the discrete adjoint method applied to structure preserving integrators have to be analyzed further as to understand the connection in a more general setting.

\section*{Acknowledgement}
This work was partly supported by the German Research Foundation (DFG, German Research Foundation) under Grant SFB 1483 – Project-ID 442419336.

\noindent
\begin{minipage}{0.8\textwidth}
This project has received funding from the European Union’s Horizon 2020 research and innovation programme under the Marie Skłodowska-Curie grant agreement No 860124. This publication reflects only the author's view and the Research Executive Agency is not responsible for any use that may be made of the information it contains.
\end{minipage}
\hfill
\begin{minipage}{0.175\textwidth}
    \centering
    \includegraphics[width=\textwidth]{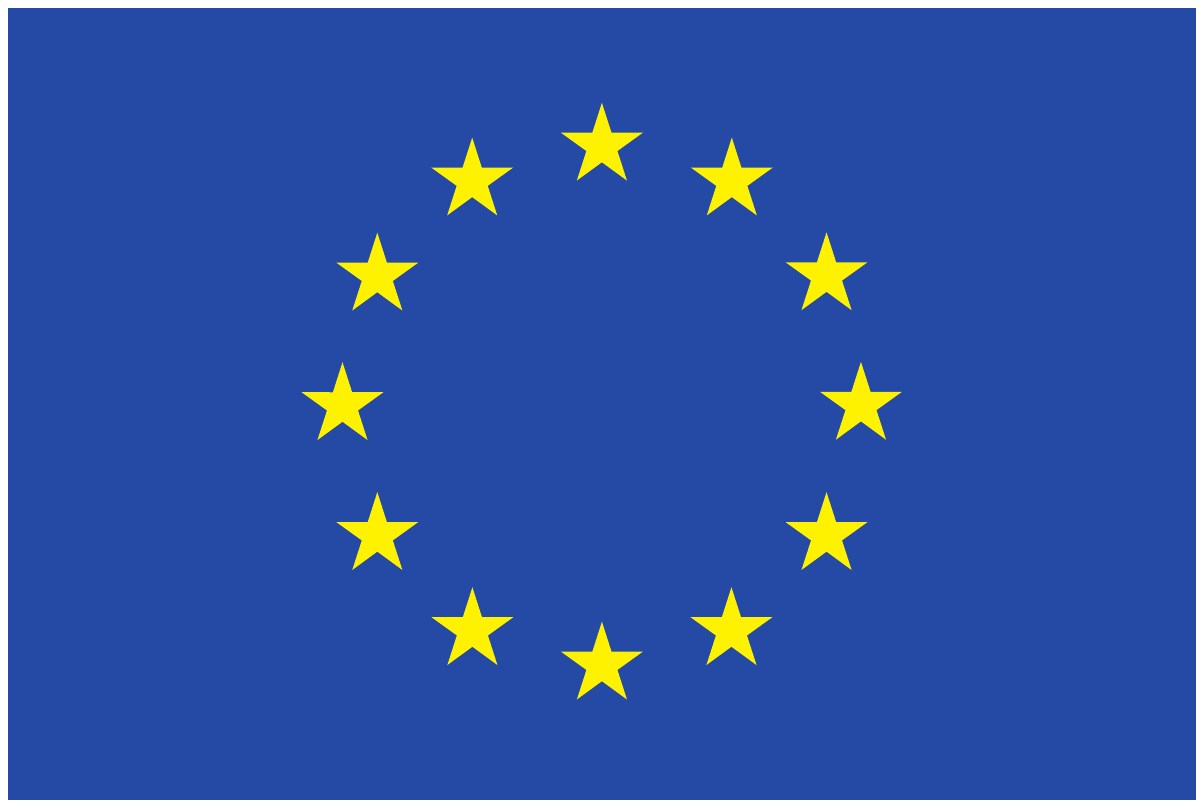}
\end{minipage}

Karin Nachbagauer acknowledges support from the Technical
University of Munich – Institute for Advanced Study.

\section{Declarations}
 
\subsection{Ethical Approval}
Not applicable\\

\subsection{Competing interests}
All authors declare that they have no conflicts of interest.\\
 
\subsection{Authors' contributions}
M.S. wrote the \textcolor{blue}{initial version of the} manuscript. R.S.T.M.A. contributed to the discussions, \textcolor{blue}{wrote much of the theoretical part of section 4 and provided additional help with figures and rewrites in other sections}. All authors reviewed the manuscript. K.N., S.O., S.L. posed the research question and conducted the first research on the topic of this paper. SL continuously supervised M.S.'s work. M.S. wrote all code.\\

\subsection{Funding}
This project has received funding from the European Union’s Horizon 2020 research and innovation programme under the Marie Skłodowska-Curie grant agreement No 860124. This publication reflects only the author's view and the Research Executive Agency is not responsible for any use that may be made of the information it contains.\\
This work was partly supported by the German Research Foundation (DFG, German Research Foundation) under Grant SFB 1483 – Project-ID 442419336.\\


\end{document}